\documentclass[hidelinks]{article}         										
\usepackage[english,strings]{babel}     
\usepackage{amsmath}   
\usepackage[utf8]{inputenc}    
\usepackage{longtable}         
\usepackage{exscale}           
\usepackage[final]{graphicx}   
\usepackage[sort]{cite}        
\usepackage{array}             
\usepackage{fancyhdr}          
\usepackage[a4paper]{geometry} 
\usepackage{xspace}            
\usepackage{tikz-cd}              
\tikzcdset{scale cd/.style={every label/.append style={scale=#1}, 
		cells={nodes={scale=#1}}}}
\usepackage{ifdraft}           
\usepackage{nchairx} 
\usepackage[sectionbib         
]{chapterbib}       
\usepackage{appendix}
\usepackage[expansion=false    
]{microtype}        
\usepackage[nottoc]{tocbibind} 
\usepackage[backref=page,      
final=true,         
pdfpagelabels       
]{hyperref}         
\usepackage{color}
\usepackage{dsfont}
\usepackage{tikz}
\usepackage{CJKutf8}
\usepackage{bbold} 
\usepackage{yfonts}
\usepackage{amssymb}

\AtEndDocument{\bigskip{\footnotesize
		\noindent\textsc{Seoul National University, Department of Mathematical Sciences, 	Research institute in Mathematics, Gwanak-Gu, 
			Seoul 08826, South Korea} \par  
		\noindent  \textit{E-mail address}: \texttt{\href{mailto:jungsoo.kang@snu.ac.kr}{jungsoo.kang@snu.ac.kr}}\par
		\noindent    \textit{E-mail address}: \texttt{\href{mailto:kevin.ruck@snu.ac.kr}{kevin.ruck@snu.ac.kr}} \par
}}

\begin{document}
	\title{Consecutive Collision Orbits in the Restricted Three-Body Problem above the First Critical Energy Value}
	\author{Jungsoo Kang and Kevin Ruck} 
	\date{}
	
	\maketitle
	
		\begin{abstract}
		In this paper, we study the planar circular restricted three-body problem for energy levels slightly above the first critical value. We first observe that the energy hypersurfaces in the Birkhoff regularization corresponding to these energy levels are of contact type. Then, using a version of Rabinowitz Floer homology, we establish the existence of either a periodic symmetric collision orbit or infinitely many symmetric consecutive collision orbits. Furthermore,  by an analytic continuation argument, for generic mass ratios and energy levels, we prove that there is no periodic symmetric collision orbit with odd number of collisions. This in turn implies the existence of at least two symmetric consecutive collision orbits. 
	\end{abstract}
	\medskip
	
	\noindent\textbf{Keywords: } Celestial Mechanics $\cdot$ Reeb Dynamics $\cdot$ Floer Homology $\cdot$ Birkhoff Regularization
	
	\medskip
	
	\noindent\textbf{Mathematics Subject Classification:} 70F16 $\cdot$ 57R17
	
	\section{Introduction}
	
	The restricted three-body problem is a classical mechanical model that describes the motion of a small body - such as a spacecraft - under the gravitational influence of two massive bodies (the primaries), which themselves move according to the solution of the two-body problem. 
	In this paper, we focus on a special case known as the planar circular restricted three-body problem. 
	This involves two simplifying assumptions: first, all motions are confined to a two-dimensional plane; and second, the two primaries move in circular orbits around their common center of mass. The first assumption is not a significant restriction, as conservation of angular momentum causes planetary motions to lie approximately within a plane. The second assumption enables us to transform the Hamiltonian from a time-dependent system into a time-independent one by adopting a rotating frame that represents the circular motion of the two primaries.
	
	\medskip

	After the pioneering work in \cite{albers2012a}, there have been many results to find special kinds of orbits by utilizing modern methods in symplectic geometry. For recent advances in the restricted three-body problem based on these methods, we refer to  \cite{frauenfelder2018a}.
	In contrast to most of the earlier works which consider energy levels below the first critical energy value, in this paper, we investigate the energy levels slightly above the first critical energy value. In this energy level, the spacecraft has enough energy to leave one of the primaries and to get to the other. To have a symplectic proof of the existence of such an orbit would be one of the ultimate goal in this research direction. However, in this paper, we will focus on consecutive collision orbits which may orbit near one of the two primaries. By consecutive collision orbits, we mean orbits that start and end in collision with one of the primaries, see \eqref{eq:collision} for the precise definition. The existence of such an orbit for energy levels below the first critical value was proved in \cite{frauenfelder2019a,ruck2024a} using Floer homology. By employing very different techniques, \cite{MORS23,ORS20} established the existence of consecutive collision orbits for a specific choice of energies. We also want to mention the related work of Chenciner and Llibre \cite{chenciner1988a} in which they proved the existence of quasi-periodic almost-collision orbits in the planar circular restricted three body problem. For a generalization of this work to the full three body problem see \cite{fejoz2002a,zhao2015a}.
	
	The planar circular restricted three-body problem carries a natural symmetry, namely the reflection about the axis connecting the two primaries. The main result of this paper is the existence of consecutive collision orbits that are symmetric with respect to this reflection. Furthermore, for a generic mass ratio and energy level, we establish the existence of at least two geometrically distinct orbits of this kind.\\[-1ex]

	\noindent\textbf{Theorem A.} {\it
		In the planar circular restricted three-body problem for energy levels slightly above the first critical energy value, the following hold. 
		\begin{enumerate}
			\item  There exists either a periodic symmetric collision orbit or infinitely many symmetric consecutive collision orbits with each primary.  
			\item For a generic choice of energy levels and mass ratios, there exist at least two geometrically distinct symmetric consecutive collision orbits with each primary.\\[-1ex]
	\end{enumerate}}

	\noindent We refer to Corollaries \ref{cor2} and \ref{cor3} for the precise formulation of the statement. 
	
	\medskip

	We outline the proof of the theorem. In Proposition \ref{contacttype}, we prove that energy hypersurfaces in the Birkhoff regularization for those energy levels are of contact type. Topologically, the regularization process corresponds to a compactification of an energy hypersurface by adding two Legendrian knots, each representing a collision with one of the primaries. Let us denote one of the Legendrian knots by $\Lambda$. Then, consecutive collision orbits correspond to Reeb chords in the Birkhoff-regularized energy hypersurface with endpoints on $\Lambda$. The regularized hypersurface inherits a symmetry  induced by the natural symmetry in the planar circular restricted three-body problem described above. A  connected component of the fixed point locus of this symmetry is also a Legendrian knot, denoted by $\Gamma$. If one can find a Reeb chord joining $\Lambda$ and $\Gamma$, then by concatenating this chord with its reflection under the symmetry, one obtains a symmetric Reeb chord with both endpoints on $\Lambda$. This Reeb chord  gives rise to a symmetric consecutive collision orbit. Therefore, our goal is to detect such Reeb chords using Rabinowitz Floer homology.

	To this end, we isotope the contact form on the Birkhoff-regularized hypersurface to a simpler one. More precisely, we show in Proposition \ref{cansubst} that the contact structure on the Birkhoff-regularized hypersurface is contactomorphic to the standard contact structure on $S^1\times S^2$. Away from the regularized locus, the regularized hypersurface is a double cover of the corresponding (unregularized) energy hypersurface. Thus, it admits a free $\mathbb{Z}_2$-action. The contactomorphism we construct is equivariant with respect to this $\mathbb{Z}_2$-action and a natural $\mathbb{Z}_2$-action on $S^1\times S^2$, and it maps the Legendrian knots $\Lambda$ and $\Gamma$ to meridians in $\{0\}\times S^2\subset S^1\times S^2$.

	Finally, in Proposition \ref{prop:equiv_rfh}, we compute the $\mathbb{Z}_2$-equivariant Lagrangian Rabinowitz Floer homology for $S^1\times S^2$ with respect to these meridians. This homology is nonzero and invariant under the contactomorphism we construct, thereby proving Theorem A.(i). In fact, this homology has infinite rank, and thus the non-existence of a symmetric periodic Reeb orbit intersecting $\Lambda$ would imply that there are infinitely many symmetric consecutive collision orbits. Note that if such a periodic Reeb orbit does exist, then its iterates could represent infinitely many homology classes. The proof of (ii) relies on an analytic continuation argument together with the fact that the rotating Kepler problem does not admit consecutive collision orbits with odd number of collisions that are also periodic after the regularization for almost all energy levels.

	\subsection*{Acknowledgments} The authors would like to thank Urs Frauenfelder and Otto van Koert for helpful discussions. This work was supported by the National Research Foundation of Korea under grants RS-2023-00211186, NRF-2020R1A5A1016126 and (MSIT) RS-2023-NR076656.

	\section{ Regularization of the Restricted Three-Body Problem}
	\label{Brik}
	For brevity, we refer to the two primaries as the Earth and the Moon, and by the restricted three-body problem (R3BP) we implicitly mean the planar circular case. 
	The Hamiltonian for the R3BP in the rotating frame is given by
	\begin{align}\label{eq:unreg_H}
		H:T^*(\mathbb{R}^2\setminus \{q^E,q^M\})\to\mathbb{R};\quad H(q,p)=\frac{1}{2}\|p\|^2+p_1q_2-p_2q_1-\frac{\mu}{\|q-q^E\|}-\frac{1-\mu}{\|q-q^M\|},
	\end{align} 
	where $q^E=(1-\mu,0)$ is the position of the Earth, $q^M=(-\mu,0)$ is the position of the Moon and $\mu$ is the mass ratio of the Earth compared to the total mass of the Earth and the Moon. It has exactly five critical points $L_1,\dots,L_5$, called the Lagrange points, ordered to satisfy $H(L_i)\leq H(L_{i+1})$. This Hamiltonian system comes with the natural involution 
	\begin{equation}
		\varrho: T^*\mathbb{R}^2 \to T^*\mathbb{R}^2;\quad  
		(q_1,q_2,p_1,p_1)\mapsto (q_1,-q_2,-p_1,p_2),
		\label{antiinvol}
	\end{equation}
	which is anti-symplectic with respect to the canonical symplectic form $\D p_1\wedge \D q_1+\D p_2\wedge \D q_2$ on $T^*\mathbb{R}^2$. Note that the fixed locus of $\varrho$ has three connected components, corresponding to
	\begin{equation}\label{eq:fix_unreg}
		\mathrm{Fix\,}\varrho \cap (\mathbb{R}^2\setminus \{q^E,q^M\})= (-\infty,q^M) \cup (q^M,q^E) \cup (q^E,\infty).	
	\end{equation}
	An orbit is said to be  symmetric if its image is invariant under $\varrho$.  
	
	An orbit $(q,p):(a,b)\to T^*(\mathbb{R}^2\setminus \{q^E,q^M\})$ of the Hamiltonian $H$ is called a {\it consecutive collision orbit} with the Earth if 
	\begin{equation}\label{eq:collision}
		\lim_{t\to a} q(t)=\lim_{t\to b} q(t)=q^E,	
	\end{equation}
	and likewise for the Moon. Note that in the literature these orbits are also sometimes referred to as \textit{ejection-collision orbits}. 
	\begin{figure}[h]
		\begin{center}
			\begin{tikzpicture}[scale=0.6,
				rocket/.pic={\draw[fill] (3.05-3.1,-0.15) .. controls (3.15-2-1.1, 0.75-1.3) .. (3.25-2-1.1, 0.95-1.3);
					\draw[fill] (3.25-2-1.1,1.35-1.3) .. controls (3.65-2-1.1, 1.25-1.3) .. (3.45-2-1.1, 1.15-1.3);
					\fill[red] (2.9-2-1.1,1.5-1.3) to[out=40, in=100, looseness=0.5] (3.4-2-1.1, 1.2-1.3) -- (3.2-2-1.1, 1-1.3) to[out=170, in=230, looseness=0.5] (2.9-2-1.1,1.5-1.3);
					\shade[
					left color = blue,
					right color = blue,
					middle color = blue!30,
					shading angle = 45
					] (0, 0) circle[radius=0.07];
					\draw (2.9-2-1.1,1.5-1.3) to[out=40, in=100, looseness=0.5] (3.4-2-1.1, 1.2-1.3);
					\draw (2.9-2-1.1,1.5-1.3) to[out=230, in=170, looseness=0.5] (3.2-2-1.1,1-1.3);
					\draw (3.4-2-1.1, 1.2-1.3) -- (3.2-2-1.1, 1-1.3);
					\draw (0, 0) circle[radius=0.07];}
				]
				\draw[arrows=->] (0,0) --(20,0);
				\draw[arrows=->] (11,-3) --(11,3);
				\draw[arrows=->] (11.5, 0.1) arc (11.3: 78.7:0.5 );
				\draw[arrows=->] (10.9, 0.5) arc (101.3: 168.7:0.5 );
				\draw[arrows=->] (10.5, -0.1) arc (191.3: 258.7:0.5 );
				\draw[arrows=->] (11.1, -0.5) arc (281.3: 348.7:0.5 );
				\path[fill] (7, 0) circle[radius=0.15];
				\path[fill] (13, 0) circle[radius=0.15];
				\draw[blue] (7, 0) to[out=120, in=90] (1, 0);
				\draw[blue] (7, 0) to[out=240, in=-90] (1, 0);
				\draw[blue] (13, 0) to[out=60, in=90] (19, 0);
				\draw[blue] (13, 0) to[out=300, in=-90] (19, 0);
				\draw (13, -2) node {\textit{Earth}};
				\draw (7, -2) node {\textit{Moon}};
				\pic[scale=1, rotate = 90] at (1.4,1) {rocket};
				\pic[scale=1, rotate = -120] at (15.05,1.45) {rocket};
			\end{tikzpicture}
		\end{center}
		\caption{Sketch of two (symmetric) consecutive collision orbits in the R3BP.}
	\end{figure}
	By regularizing collisions, consecutive collision orbits can be interpreted as orbits that start and end in the cotangent fiber over the collision point. If a consecutive collision orbit closes up after regularization, we refer to it as a {\it periodic collision orbit}. 
	To make this precise, we need to regularize the Hamiltonian system at an energy level $c\in\mathbb{R}$ of our interest, meaning that we compactify the noncompact part of $H^{-1}(c)$ arising from  collisions and extend the dynamics accordingly. The first step is to reparametrize the time variable, for example
	\begin{align*}
		t_{\text{new}}=\int\frac{1}{\| q-q^E\|}\D t,
	\end{align*}
	which, in terms of the Hamiltonian function, corresponds to redefining it as
	\begin{align*}
		H_{\text{new}}=\| q-q^E\|(H-c).
	\end{align*}
	It is easy to see that this will regularize the divergence at $q=q^E$. The next step in this procedure now depends on a specific regularization one wishes to employ. In this section, we first recall the Moser regularization \cite{moser1970a}, which regularizes one of the collision points, and review the foundational results in \cite{albers2012a}. Then, we adapt these results to the setting of the Birkhoff regularization, which regularizes both the Earth and the Moon collisions simultaneously. 
	
	\subsection{Moser regularization}
	\label{basics}

	Let $c<H(L_1)$ be an energy level below the first critical value of $H$. Then, the energy hypersurface $H^{-1}(c)$ consists of three connected components: one around the Earth $\Sigma_c^E$, one around the Moon $\Sigma_c^M$, and one far from both bodies. In performing a regularization, we compactify only the components $\Sigma_c^E$ and $\Sigma_c^M$, while disregarding the third component, as for all relevant energy levels, a spacecraft starting near the Earth or the Moon cannot access this region of space. To regularize collisions at the Earth, we first translate the coordinates so that the Earth is located at the origin. We denote the Earth component again by $\Sigma_c^E$. Then, we switch the base and fiber coordinates of the cotangent bundle $T^*\mathbb{R}^2$ via the map
	\begin{equation}
		\mathfrak{sw}: T^*\mathbb{R}^2 \to T^*\mathbb{R}^2 \  ;\quad (q_1,q_2,p_1,p_2)\mapsto (-p_1,-p_2,q_1,q_2).
	\end{equation}  
	After applying the switch map, we map $\mathfrak{sw}(\Sigma_c^E)$ to the cotangent bundle of the two-sphere via the cotangent lift $T^*\phi$ of the inverse of the stereographic projection:
	\[
	\mathfrak{\phi}: \mathbb{R}^2 \to S^2 \subset \mathbb{R}^3; \quad x=(x_1,x_2)\mapsto \left(\frac{2x_1}{1+\|x\|^2}, \frac{2x_2}{1+\|x\|^2}, \frac{\|x\|^2-1}{1+\|x\|^2}\right).
	\] 
	We abbreviate the Morse regularization map by 
	\[
	\mathfrak{M}:=T^*\phi\circ\mathfrak{sw}:T^*\mathbb{R}^2\to T^*S^2.
	\]
	Since both $T^*\phi$ and $\mathfrak{sw}$ are symplectic, so is $\mathfrak{M}$. The compactification of $\mathfrak{M}(\Sigma_c^E)\subset T^* S^2$, denoted by $\overline{\Sigma}_c^{\mathfrak{M},E}$, is the desired regularized hypersurface. The non-compactness of the energy hypersurface $\Sigma_c^E$ stems from the momentum going to infinity as one approaches the collision point $q^E$. Through the procedure described above, we map the fibers (momenta) of $T^*\mathbb{R}^2$ to the base (positions) of $T^*S^2$, where the infinite momentum corresponds to the north pole of $S^2$. The Moon component $\Sigma_c^M$ can be regularized in the same manner, and we write $\overline{\Sigma}_c^{\mathfrak{M},M}$ for the regularized hypersurface.


	\begin{theorem}(\!\!{\cite{albers2012a}})\label{thm:moser}
		For every $c<H(L_1)$, the Moser-regularized hypersurface $\overline{\Sigma}_c^{\mathfrak{M},E}$ is transverse to the standard Liouville vector field on $T^*S^2$. In particular, the restriction of the Liouville one-form on $T^*S^2$ to $\overline{\Sigma}_c^{\mathfrak{M},E}$ is a contact form, and the associated Reeb flow  corresponds to a reparametrization of the Hamiltonian flow of the R3BP. The same statement holds  for $\overline{\Sigma}_c^{\mathfrak{M},M}$.
	\end{theorem}
	
	Topologically, the compactification is done by adding a circle fiber over the north pole $n\in S^2$. We denote this fiber by 
	\[
	\Lambda_{c}^{\mathfrak{M},E} = \overline{\Sigma}_c^{\mathfrak{M},E} \cap T^*_nS^2,
	\] 
	which is a Legendrian knot. 
	Then, consecutive collision orbits with the Earth at energy level $c$ correspond to Reeb chords on $\overline{\Sigma}_c^{\mathfrak{M},E}$ that begin and end on $\Lambda_{c}^{\mathfrak{M},E}$. We write $\Lambda_{c}^{\mathfrak{M},M}$ for the circle fiber in the case of the Moon. Note that according to \cite[Lemma~5.7.1]{frauenfelder2018a} these fibers are convex.
	
	Through the Moser regularization, the anti-symplectic involution in \eqref{antiinvol} corresponds to the anti-symplectic involution 
	\[
	\mathfrak{i}\circ T^*\mathfrak{r}:T^*S^2\to T^*S^2,
	\] 
	where $\mathfrak{i}$ maps $(x,\zeta)\in T^*S^2$ to $(x,-\zeta)$ and $T^*\mathfrak{r}$ is the cotangent lift of the reflection $\mathfrak{r}$ on $S^2$ about the meridian $S^2\cap(\{0\}\times\mathbb{R}^2)$. Note that the fixed locus $\mathrm{Fix\,}(\mathfrak{i}\circ T^*\mathfrak{r})$ is 
	the conormal bundle $N_{\mathrm{Fix\,}\mathfrak{r}}^*S^2$ over the meridian  $\mathrm{Fix\,}\mathfrak{r}$.
	The fixed locus $\mathrm{Fix\,}\varrho\cap \Sigma_c^E$ corresponds to the fixed locus $\mathrm{Fix\,}(\mathfrak{i}\circ T^*\mathfrak{r})\cap \overline{\Sigma}_c^{\mathfrak{M},E}$ via $\mathfrak{M}$. The latter is the union of two circles since $\overline{\Sigma}_c^{\mathfrak{M},E}$ is fiberwise starshaped by Theorem \ref{thm:moser}. Each circle corresponds to the component of $\mathrm{Fix\,}\varrho\cap \Sigma_c$ with $q_1$ on the left- or the right-hand side of the Earth, see \eqref{eq:fix_unreg}. We take the component corresponding to the right-hand side of the Earth, which is further from the Moon, and denote it by
	\[
	\Gamma_c^{\mathfrak{M},E}\subset \overline{\Sigma}_c^{\mathfrak{M},E} \cap \mathrm{Fix\,}(\mathfrak{i}\circ T^*\mathfrak{r}).
	\]
	It is also a Legendrian knot.  
	Moreover, $\Lambda_{c}^{\mathfrak{M},E}$ and $\Gamma_c^{\mathfrak{M},E}$ intersect exactly once.

	\medskip
	
	The energy level slightly above the first critical value is also studied in  \cite{albers2012a}. As the energy level surpasses $H(L_1)$, the Earth and the Moon components become connected.  At the critical level $H(L_1)$, the Lagrange point $L_1$ divides this bounded component of $H^{-1}(L_1)$ into two connected components. We denote by $\Sigma_{H(L_1)}^E \subset H^{-1}(L_1)$ the union of the Earth component and $L_1$. Again, $\mathfrak{M}(\Sigma_{H(L_1)}^E)\cap \mathrm{Fix\,}(\mathfrak{i}\circ T^*\mathfrak{r})$ has two connected components, and we write $\widehat \Gamma_{H(L_1)}^{\mathfrak{M},E}$ for the component left-hand side of the Earth, which is closer to the Moon. Then, we have
	\begin{equation}\label{eq:lagrange}
		\mathfrak{M}(L_1)=T_s^*S^2\cap \widehat \Gamma_{H(L_1)}^{\mathfrak{M},E},	
	\end{equation}
	where $T_s^*S^2$ denotes the cotangent fiber at the south pole $s\in S^2$.
	
	For $d\in (H(L_1),H(L_2))$,  we denote the bounded component of $H^{-1}(d)$  by $\Sigma_d$. 
	
	\begin{theorem}(\!\!{\cite{albers2012a}})\label{thm:moser2}
		There exists $\epsilon>0$ such that, for every $d\in(H(L_1),H(L_1)+\epsilon)$, the Moser-regularized hypersurface $\overline{\Sigma}_d^\mathfrak{M}$ admits a contact form whose Reeb flow  corresponds to a reparametrization of the Hamiltonian flow of the R3BP.
	\end{theorem}
	
	Let $d\in(H(L_1),H(L_1)+\epsilon)$ be as in the above theorem. To regularize $\Sigma_d$, we divide it into three pieces, the neck region (the intersection of $\Sigma_d$ and an open ball of the Lagrange point $L_1$ in $T^*(\mathbb{R}^2\setminus\{q^E,q^M\})$), the Earth component, and the Moon component. Then, we compactify the Earth and the Moon components as in the case of energy level $c<H(L_1)$ using the Moser regularization. Note that the resulting regularized hypersurface $\overline{\Sigma}_d^\mathfrak{M}$ is topologically $\Sigma_d$ compactified by adding two circles.

	To show the contact property of $\overline{\Sigma}_d^\mathfrak{M}$, the authors of \cite{albers2012a} observe that the Earth and the Moon components (after the Moser regularization) are still transverse to the standard Liouville vector field on $T^* S^2$ as in Theorem \ref{thm:moser}. Furthermore, they proved that this Liouville vector field can be extended over the neck part. Thus, $\overline{\Sigma}_d^\mathfrak{M}$ can be viewed as  a contact connected sum of $\overline{\Sigma}_{c}^{\mathfrak{M},E}$ and $\overline{\Sigma}_{c}^{\mathfrak{M},M}$ for $c<H(L_1)$, where the connected sum operation produces the neck region. 
	
	We denote the Legendian knots in $\overline{\Sigma}_d^\mathfrak{M}$ that correspond to $\Lambda_{c}^{\mathfrak{M},E},\Gamma_c^{\mathfrak{M},E}\subset \overline{\Sigma}_c^{\mathfrak{M},E}$ by
	\[
	\Lambda_{d}^{\mathfrak{M},E},\;\Gamma_{d}^{\mathfrak{M},E} \subset \overline{\Sigma}_d^\mathfrak{M},
	\]
	respectively. 
	More precisely, $\Lambda_{d}^{\mathfrak{M},E}$ is the intersection of $T^*_nS^2$ with the Moser-regularized Earth part of $\Sigma_d$. The involution $\varrho$ restricted to ${\Sigma_d}$ extends to $\overline{\Sigma}  _d^\mathfrak{M}$. The fixed locus of this extended involution has three connected components, that are characterized by the position of $q_1$, see \eqref{eq:fix_unreg}. Then, $\Gamma_{d}^{\mathfrak{M},E}$ is the component corresponding to the case of $q_1>q^E$. We remark that, due to \eqref{eq:lagrange}, the Legendrians knots $\Lambda_{d}^{\mathfrak{M},E}$ and $\Gamma_{d}^{\mathfrak{M},E}$ are away from the neck region and thus indeed subsets of $\overline{\Sigma}_d^\mathfrak{M}$.

	The standard Liouville flow on $T^*S^2$ gives a contactomorphism between each of $\overline{\Sigma}_{c}^{\mathfrak{M},E}$ and $\overline{\Sigma}_{c}^{\mathfrak{M},M}$ to the unit cotangent bundle $U^*S^2$ for $c<H(L_1)$. Furthermore, the contactomorphism $\overline{\Sigma}_{c}^{\mathfrak{M},E}\to U^*S^2$ maps $\Lambda_{c}^{\mathfrak{M},E}$ to $U^*_nS^2$, the fiber over the north pole $n\in S^2$. The Legendrian knot  $\Gamma_{c}^{\mathfrak{M},E}$ is mapped to one  component of the unit normal bundle over the meridian $S^2\cap(\{0\}\times\mathbb{R}^2)$. If we view $T^*S^2$ as a submanifold of $\mathbb{R}^3\times(\mathbb{R}^3)^*\cong \mathbb{R}^3\times\mathbb{R}^3$, this component is written as $\Gamma:=(S^2\cap(\{0\}\times\mathbb{R}^2))\times \{(1,0,0)\}$. For $d\in (H(L_1),H(L_1)+\epsilon)$, we also have a contactomorphism 
	\begin{equation}
		\begin{aligned}
			\overline{\Sigma}_d^\mathfrak{M}  \longrightarrow U^*S^2\# U^*S^2,
		\end{aligned}
		\label{eq:cont_conn}
	\end{equation}
	where $\#$ refers to the contact connected sum, which maps $\Lambda_{d}^{\mathfrak{M},E}$ and $\Gamma_{d}^{\mathfrak{M},E}$ to $U_n^*S^2$ and one component of the unit normal bundle over the meridian $S^2\cap(\{0\}\times\mathbb{R}^2)$ in the first $U^*S^2$, respectively. Such Legendrian knots in $U^*S^2\# U^*S^2$ are away from the region where the connected sum operation is performed.

	\subsection{Birkhoff Regularization}
	\label{IntroBirk}
	In this section, we introduce the Birkhoff regularization \cite{birkhoff1915a}. See also \cite[Section 4]{CFZ23}   for an insightful discussion on this. In contrast to the Moser regularization,  the Birkhoff regularization regularizes collisions with both primaries simultaneously. It is given by the cotangent lift
	\[
	\mathfrak{B}:=T^*\mathfrak{b}: T^*\left(\mathbb{C}\setminus \left\{-\tfrac{1}{2}, 0, \tfrac{1}{2}\right\}\right) \to T^*\mathbb{C}; \ \ (z,w)\mapsto \left(\frac{z^2+\frac{1}{4}}{2z}, \frac{2\bar{z}^2w}{\bar{z}^2-\frac{1}{4}}\right)
	\label{Bmap}
	\]
	of the map
	\[
	\mathfrak{b}: \mathbb{C}\setminus \{0\} \to \mathbb{C}; \ \ b(z)= \frac{1}{4}\left(2z+\frac{1}{2z}\right).
	\label{bmap}
	\]
	It regularizes the restricted three-body problem with primaries at $(\pm \frac{1}{2},0)$ by pulling back the Hamiltonian $H$ in \eqref{eq:unreg_H} shifted such that the Earth and the Moon are located at $(\frac{1}{2},0)$ and $\left(-\frac{1}{2},0\right)$, respectively. We write $H$ again for this shifted Hamiltonian. More precisely, the regularized Hamiltonian for the energy level $d\in\mathbb{R}$ is given by 
	\[
	K: T^*\left(\mathbb{C}\setminus\{0\}\right) \to \mathbb{R}\ ;\quad 
	K(z,w)=  \mathfrak{B}^*\left(\Big|q-\frac{1}{2}\Big|\Big|q+\frac{1}{2}\Big|\big(H(q,p)-d\big)\right)(z,w).
	\]
	It is explicitly written as follows:
	\[
	\begin{aligned}
		K(z,w)= &\; \frac{|z|^2|w|^2}{2} - \frac{\mu \left|z-\frac{1}{2}\right|^2}{2|z|} - \frac{(1-\mu) \left|z+\frac{1}{2}\right|^2}{2|z|}\\ &\quad+\frac{\left|z-\frac{1}{2}\right|^2\left|z+\frac{1}{2}\right|^2}{4|z|^2}\left(\text{Im}\left(z\bar{w}\frac{z+\frac{1}{2}}{z-\frac{1}{2}}-\frac{2\mu z^2\bar{w}}{z^2-\frac{1}{4}}\right)-d\right).    
	\end{aligned}
	\]
	We remark that the Hamiltonian $K$ is not singular at $z=\pm\frac{1}{2}$. Moreover, the singular point $z=0$ corresponds to $q=\infty$, and therefore does not interfere with the energy hypersurface we consider. We are particularly interested in the case of $d\in(H(L_1),H(L_2))$, which we assume from now on. The bounded component of $K^{-1}(0)$, denoted by $\overline{\Sigma}_d^{\mathfrak B}$, is the Birkhoff-regularized energy hypersurface of the bounded component $\Sigma_d$ of $H^{-1}(d)$. 	The Birkhoff regularization adds the following two fiber circles over the Earth and the Moon to $H^{-1}(d)$,
	\[
	\Lambda_d^{\mathfrak{B},E}:=\overline{\Sigma}_d^\mathfrak{B} \cap T^*_{\frac{1}{2}}\mathbb{C},\qquad	\Lambda_d^{\mathfrak{B},M}:=\overline{\Sigma}_d^\mathfrak{B} \cap T^*_{-\frac{1}{2}}\mathbb{C}.
	\]
	
	The map $\mathfrak{b}$ is a branched double covering map with the branch locus $z=\pm\frac{1}{2}$. This leads to an additional $\mathbb{Z}_2$-symmetry in the Hamiltonian system of $K$. More precisely, the action of $\mathbb{Z}_2=\{-1,+1\}$ given by 
	\[
	1\acts (z,w):= (z,w)\ \ \text{ and }\ \ -1\acts (z,w):= \left(\frac{1}{4z}, -4\bar{z}^2w\right), 
	\]
	preserves the standard symplectic form on $T^*\mathbb{C}$, and the map $\mathfrak{B}$ is invariant under this action, i.e.~$\mathfrak{B}(z,w)=\mathfrak{B}\left(-1\acts (z, w)\right)$. This implies that the regularized energy hypersurface $\overline{\Sigma}^{\mathfrak B}_d$ admits an involution under which the dynamics remain invariant. As this involution arises from the fact that the Birkhoff regularization relies on the branched double covering $\mathfrak{b}$, we will quotient it out when considering Rabinowitz Floer homology for $\overline{\Sigma}^{\mathfrak B}_d$.

	As observed in  \cite[Proposition 4.3]{CFZ23}, $\overline{\Sigma}_d^{\mathfrak{B}} $ is diffeomorphic to $S^1\times S^2$, and thus is a double cover of $\overline{\Sigma}_d^{\mathfrak{M}} $, which is diffeomorphic to $\mathbb{R}P^3\# \mathbb{R}P^3$, see \eqref{eq:cont_conn}. In the next section, we observe that this double covering map can be arranged to respect the dynamics on both hypersurfaces.

	The natural involution $\varrho$ in \eqref{antiinvol} lifts to two distinct anti-symplectic involutions:
	\begin{align}
		\varrho_1(z,w):= (\bar{z}, -\bar{w}),\qquad \varrho_2(z,w):= \left(\frac{1}{4\bar{z}}, 4z^2\bar{w}\right).
	\end{align}
	The fixed locus of $\varrho_1$ is $\mathbb{R}\times i\mathbb{R} \subset \mathbb{C}^2=T^*
	\mathbb{C}$, corresponding to the components $(-\infty,q^M)\cup (q^E,\infty)$ in \eqref{eq:fix_unreg}. The fixed locus of $\varrho_2$ is a double twisted M\"obius strip, corresponding to the component $(q^M,q^E)$. The regularized hypersurface $\overline{\Sigma}^{\mathfrak B}_d$ and the fiber circles $\Lambda_d^{\mathfrak{B},E}$, $\Lambda_d^{\mathfrak{B},M}$ are invariant under both involutions. We are only interested in $\varrho_1$ and denote 
	\[
	\mathrm{Fix\,}\varrho_1\cap \overline{\Sigma}^{\mathfrak B}_d =  \Gamma_d^{\mathfrak{B},M}\cup \Gamma_d^{\mathfrak{B},E},
	\]
	where $\Gamma_d^{\mathfrak{B},M}$ and $\Gamma_d^{\mathfrak{B},E}$ are the connected components corresponding to $(-\infty,q^M)$ and $(q^E,\infty)$, respectively. 
	
	\begin{figure}[h]
		\begin{center}
			\begin{tikzpicture}
				\draw (-6,0) --(3,0);
				\draw[->] (0,-2) -- (0,2);
				\draw[purple!80!black][ultra thick] (-6,0)--(-3,0);
				\draw[blue!80!black][ultra thick] (0,0)--(3,0);
				\draw[gray!80!black] (-3, 0) node{$\bullet$};
				\draw[gray!80!black] (-3, 0) circle[radius=0.2];
				\fill[gray!80!black,opacity=0.3] (-3, 0) circle[radius=0.2];
				\draw[gray!80!black] (0, 0) node{$\bullet$};
				\draw[green!60!black] (0, 0) circle[radius=0.2];
				\fill[green!60!black,opacity=0.3] (0, 0) circle[radius=0.2];
				\draw[purple!80!black] (-3.9,-1.2) node{$\leftrightarrow \Gamma_d^{\mathfrak{B},M}$};
				\draw[purple!80!black] (-6,-0.3) to[out=270, in=90, looseness=1] (-4.5,-1);
				\draw[purple!80!black] (-3,-0.3) to[out=270, in=90, looseness=1] (-4.5,-1);
				\draw[blue!80!black] (3,-0.3) to[out=270, in=90, looseness=1] (1.5,-1);
				\draw[blue!80!black] (0,-0.2) to[out=270, in=90, looseness=1] (1.5,-1);
				\draw[blue!80!black] (2,-1.2) node{$\leftrightarrow\Gamma_d^{\mathfrak{B},E}$};
				\draw[gray!50!black] (-2.5,0.8) node{Moon $\leftrightarrow \Lambda_d^{\mathfrak{B},M}$};
				\draw[green!50!black] (1.5,0.8) node{Earth $\leftrightarrow \Lambda_d^{\mathfrak{B},E}$};
			\end{tikzpicture}
		\end{center}
		\caption{Interpretation of the above define Legendrian submanifolds in the position space picture of the unregularized setting.}
	\end{figure}
	\begin{remark}\label{rem:double}
		Suppose that there is a Hamiltonian chord $\gamma:[0,T]\to \overline{\Sigma}^{\mathfrak B}_d$ with $\gamma(0) \in \Lambda_d^{\mathfrak{B},E}$ and $\gamma(T)\in \Gamma_d^{\mathfrak{B},E}$. Then, the curve  $t\mapsto \varrho_1\circ\gamma(T-t)$ is also a Hamiltonian chord with $\varrho_1\circ\gamma(T)=\gamma(T)$ and $\varrho_1\circ\gamma(0)\in \varrho_1(\Lambda_d^{\mathfrak{B},E})=\Lambda_d^{\mathfrak{B},E}$. Concatenating this with $\gamma$, we obtain a $\varrho_1$-symmetric chord with both endpoints in $\Lambda_d^{\mathfrak{B},E}$. This corresponds to a symmetric consecutive collision orbit with the Earth in the (unregularized) R3BP.
		
		Similarly, if there is a Hamiltonian chord $\gamma:[0,T]\to \overline{\Sigma}^{\mathfrak B}_d$ with both endpoints in $\Gamma_d^{\mathfrak{B},E}$, then by concatenating it with $t\mapsto \varrho_1\circ\gamma(T-t)$, we obtain a $\varrho_1$-symmetric periodic orbit, which corresponds to a symmetric periodic (possibly collision) orbit near the Earth in the R3BP.
	\end{remark}

	\subsection{Contact property of Birkhoff regularized energy hypersurfaces}
	\label{secContact}
	In order to apply Rabinowitz Floer homology, a tool from contact and symplectic geometry, we want to prove that the Birkhoff-regularized energy hypersurfaces $\overline{\Sigma}^\mathfrak{B}_d$ are of contact type. A corresponding result is proved in \cite{albers2012a} for the Moser-regularized energy hypersurfaces as recalled in Theorems \ref{thm:moser} and \ref{thm:moser2}. To carry this result over to our setting, we construct a map relating the Birkhoff regularization and the Moser regularization.

	To have an explicit expression of $T^*\phi$, we first compute the differential of $\phi$:
	\begin{align}
		\D \phi_x= \frac{1}{(\|x\|^2+1)^2} \left(
		\begin{matrix*}
			2\|x\|^2-4x_1^2+2 & -4x_1x_2\\
			-4x_1x_2 & 2\|x\|^2-4x_2^2+2\\
			4x_1 & 4x_2
		\end{matrix*}\right) : T_{x}\mathbb{R}^2 \to T_{\phi(x)}S^2.
	\end{align}
	We then trivialize the image of $d\phi$ (away from the south pole) using the following orthogonal transformation:
	\[
	A_x:= \frac{1}{\|x\|^4+\|x\|^2}\left(
	\begin{matrix*}
		2x_2^2+\|x\|^4-\|x\|^2 & -2x_1x_2 & -2\|x\|^2 x_1\\
		-2x_1x_2 & 2x_1^2+\|x\|^4-\|x\|^2 & -2\|x\|^2 x_2\\
		2\|x\|^2 x_1 & 2\|x\|^2 x_2 & \|x\|^4-\|x\|^2
	\end{matrix*}\right)
	\]
	One can think of this map as the rotation of a plain tangent to $S^2$ at the point $\phi(x)$ into $\mathbb{R}^2\subset\mathbb{R}^3$. Indeed,  a direct computation shows 
	\[
	A_x\circ d\phi_x = \frac{2}{\|x\|^4+\|x\|^2} \left(\begin{matrix}
		x_2^2-x_1^2 & -2x_1x_2\\
		-2x_1x_2 & x_1^2-x_2^2\\
		0 & 0\\
	\end{matrix}\right).
	\]
	To simplify computations, we introduce the complex coordinates $x=x_1+ix_2$ and $v=v_1+iv_2$. Then, 
	\[
	A_x\circ d\phi_x(v)= \frac{-2}{\|x\|^4+\|x\|^2}x^2\bar{v}
	\]
	Remembering that $T^*\phi(x,\xi)=(\phi(x),((d\phi_x)^{-1})^*\xi )$, we apply the inverse and the complex conjugate to the trivialized differential $A_x\circ d\phi_x$ and obtain 
	\[
	\xi \mapsto -\frac{\|x\|^2+1}{2\|x\|^2}x^2\bar{\xi }.
	\]
	Therefore, the cotangent lift $T^*\phi$ with respect to the trivialization $A_x$ is written as 
	\[
	(x,\xi)\mapsto \Bigg(\phi(x), -\frac{\|x\|^2+1}{2\|x\|^2}x^2\bar{\xi }\Bigg).
	\]

	\medskip
	
	Next, we want to see that a suitable composition of the Moser regularization map $\mathfrak{M}=T^*\phi\circ\mathfrak{sw}$ and  the Birkhoff regularization map $\mathfrak{B}$ extends over the collision loci $\pm\frac{1}{2}$. To simplify computations, we shift to the coordinate $z= z^\prime-\frac{1}{2}$ when analyzing the collision at $-\frac{1}{2}$ (analogously for the collision at $\frac{1}{2}$) and study the collision at $z^\prime= 0$ instead. In addition, since the Moser regularization requires the collision locus to be at the origin, we shift the image of $\mathfrak{B}$ by $\frac{1}{2}$ ($-\frac{1}{2}$ for the collision at $\frac{1}{2}$). In summary, the modified the Birkhoff regularization map is given by
	\[
	\mathfrak{B}':(z^\prime, w)\mapsto \left(\frac{{z^\prime}^2}{2z^\prime-1}, \frac{(4\bar{z^\prime}^2-4\bar{z^\prime}+1)w}{2\bar{z^\prime}^2-2\bar{z^\prime}}\right).
	\]
	For ease of notation, we denote $z'$ simply by $z$. We compute
	\begin{equation*}
		\begin{aligned}
			\mathfrak{M}\circ \mathfrak{B}'  (z,w)&= T^*\mathfrak{\phi}\left(\frac{(4\bar{z}^2-4\bar{z}+1)w}{2\bar{z}^2-2\bar{z}},\, \frac{z^2}{2z-1}\right)\\
			&=\left(\phi\Big(\frac{(4\bar{z}^2-4\bar{z}+1)w}{2\bar{z}^2-2\bar{z}}\Big),\,  -\frac{1+\frac{\left|2\bar{z}^2-2\bar{z}\right|^2}{\left|4\bar{z}^2-4\bar{z}+1\right|^2|w|^2}}{2} 
			\frac{\left(4\bar{z}^2-4\bar{z}+1\right)^2 w^2}{\left(2\bar{z}^2-2\bar{z}\right)^2} \frac{\bar{z}^2}{2\bar{z}-1}\right).
		\end{aligned}
	\end{equation*}
	A straightforward computation shows that this map extends smoothly and regularly over $z=0$. Moreover, 
	\[
	\mathfrak{M}\circ \mathfrak{B'}(0,w)=\Big(n, \frac{w^2}{8}\Big),
	\]
	where $n\in S^2$ is the north pole, which shows that $\Lambda_{d}^{E,\mathfrak{B}}$ is a double cover of $ \Lambda_{d}^{E,\mathfrak{M}}$ via this map.  An analogous computation holds for the other collision point with $\frac{w^2}{8}$ replaced by $-\frac{w^2}{8}$.

	From the above observation, we can conclude that there is a natural double covering map 
	\begin{equation}\label{eq:cover_moser}
		\pi^{\mathfrak{B},\mathfrak{M}}:\overline{\Sigma}_d^{\mathfrak{B}} \longrightarrow \overline{\Sigma}_d^{\mathfrak{M}}.		
	\end{equation}
	Indeed, we have the composition of two smooth maps
	\[
	{\Sigma}_d^{\mathfrak{B}}=\overline{\Sigma}_d^{\mathfrak{B}}\setminus(\Lambda_d^{\mathfrak{B},E}\cup \Lambda_d^{\mathfrak{B},M})\longrightarrow \Sigma_d \longrightarrow \overline{\Sigma}_d^{\mathfrak{M}}\setminus(\Lambda_d^{\mathfrak{M},E}\cup\Lambda_d^{\mathfrak{M},M}),
	\]	
	where the first map is the double covering map in the Birkhoff regularization and the second one is a diffeomorphism in the Moser regularization. Moreover, we observed above that this composition extends over the collision loci, and gives rise to a local diffeomorphism $\pi^{\mathfrak{B},\mathfrak{M}}$, which is hence the claimed double covering map.

	By lifting the contact form on $\overline{\Sigma}_d^\mathfrak{M}$ in Theorem \ref{thm:moser2} to $\overline{\Sigma}_d^\mathfrak{B}$ via $\pi^{\mathfrak{B},\mathfrak{M}}$, we prove the following proposition.
	
	\begin{proposition}\label{prop:contact_birkhoff}
		Let $\epsilon>0$ be as in Theorem \ref{thm:moser2}. Then for every $d\in(H(L_1),H(L_1)+\epsilon)$, the Birkhoff-regularized energy hypersurface $\overline{\Sigma}_d^\mathfrak{B}$ admits a contact form whose Reeb flow corresponds to a reparametrization of the Hamiltonian flow of the R3BP.
		\label{contacttype}
	\end{proposition}

	\section{Reeb dynamics on $S^1\times S^2$}
	\label{3}
	
	The goal of this section is twofold. We first construct a contactomorphism between  the Birkhoff-regularized hypersurface $\overline{\Sigma}_d^\mathfrak{B}$ for $d\in(H(L_1),H(L_1)+\epsilon)$ and $S^1\times S^2$ with the standard contact structure, which is equivariant under  natural symplectic involutions in both spaces. This contactomorphism maps Legendian knots $\Lambda_d^{\mathfrak{B},E}$ and $\Lambda_d^{\mathfrak{B},M}$, which are of interest to us, to meridians in $\{0\}\times S^2\subset S^1\times S^2$ and $\{\pi\}\times S^2\subset S^1\times S^2$, respectively. Then, we compute the Robbin-Salamon indices of Reeb chords with endpoints in those meridians for the standard Reeb flow on $S^1\times S^2$. These results will be used in the next section when computing a relevant Rabinowitz Floer homology.

	\subsection{Contactomorphism to the standard contact structure}
	\label{3.1}
	
	Recall from above \eqref{eq:cont_conn} that $U^*_nS^2$ is the unit cotangent fiber over the north pole $n\in S^2$ and $\Gamma:=(S^2\cap(\{0\}\times\mathbb{R}^2))\times \{(1,0,0)\}$.
	\begin{lemma}\label{lem:covering}
		There is a double covering map from $\pi:S^3\to U^*S^2$ such that $\pi^{-1}(U^*_nS^2)$ and $\pi^{-1}(\Gamma)$ are contained in the two-sphere $S^3\cap \mathrm{span}\{1,\mathbf{i},\mathbf{j}\}$, where we identified $\mathbb{R}^4$ with the vector space of quaternions $\mathbb{H}=\mathrm{span}\{1,\mathbf{i},\mathbf{j},\mathbf{k}\}$. Furthermore, the pullback of the contact structure of $U^*S^2$ by $\pi$ is the standard contact structure on $S^3$.
	\end{lemma}
	\begin{proof}
		
		Through the diffeomorphism 
		\[
		\Xi:U^*S^2 \overset{\cong}{\longrightarrow} SO(3);\quad (p,v)\mapsto(v\ p\ v\times p),
		\]
		the circles $\Gamma$ and $U_n^*S^2$ are mapped as follows:
		\begin{align*}
			\Gamma \mapsto 
			\left\{\begin{pmatrix}
				1&0&0\\
				0&\cos(\theta)&-\sin(\theta)\\
				0&\sin(\theta)&\cos(\theta)
			\end{pmatrix}\ \bigg\vert\ \theta\in [0,2\pi)\right\} \subset SO(3)
		\end{align*}
		and
		\begin{align*}
			U_n^*S^2 \mapsto 
			\left\{\begin{pmatrix}
				\cos(\theta)&0&\sin(\theta)\\
				\sin(\theta)&0&-\cos(\theta)\\
				0&1&0
			\end{pmatrix}\ \bigg\vert\ \theta\in [0,2\pi)\right\} \subset SO(3).
		\end{align*}

		To construct a double covering map, we interpret $S^3$ as the group of unit quaternions. Recall that a unit quaternion represents a rotation in $\mathbb{R}^3$ like $SO(3)$. We identify the imaginary part of the vector space of quaternions with $\mathbb{R}^3$, i.e.~$\mathbb{R}^3 \cong \text{Im}(\mathbb{H})$. The rotation associated to the rotation axis $u\in\text{Im}(\mathbb{H})$ and the rotation angle $\vartheta\in[0,2\pi)$ is given by 
		\[
		R_q:\text{Im}(\mathbb{H}) \to \text{Im}(\mathbb{H});\quad  {v}\mapsto q{v}q^{-1}
		\]
		where
		\begin{align*}
			q:=\cos\left(\frac{\vartheta}{2}\right)+\sin\left(\frac{\vartheta}{2}\right)u \in \mathbb{H}.
		\end{align*}
		This interpretation gives the double covering map $\tilde\pi:S^3\to SO(3)$, which sends $q$ and $-q$ to $R_q\in SO(3)$. We observe 
		\[
		\tilde{\pi}\big(\cos(\tfrac{\theta}{2})+\sin(\tfrac{\theta}{2})\mathbf{i}\big)=\begin{pmatrix}
			1&0&0\\
			0&\cos(\theta)&-\sin(\theta)\\
			0&\sin(\theta)&\cos(\theta)
		\end{pmatrix}\]
		whose rotation axis is spanned by $(1,0,0)$. 
		The matrix $\begin{pmatrix}
			\cos(\theta)&0&\sin(\theta)\\
			\sin(\theta)&0&-\cos(\theta)\\
			0&1&0
		\end{pmatrix}$ in $\Xi(U_n^*S^2)$ is the rotation by $\arccos(\frac{\cos(\theta)-1}{2})$ with respect to the rotation axis spanned by $(\frac{1+\cos(\theta)}{\sin(\theta)},1,1)$. The rotation axis of $R_q$ is spanned by $\mathrm{Im}(q)=\frac{1}{2}(q-\bar q)$. Therefore, 
		\begin{equation}
			\tilde{\pi}^{-1}(\Xi(\Gamma\cup U^*_nS^2))\subset S^3\cap \mathrm{span}\{1,\mathbf{i},\mathbf{j}+\mathbf{k}\}.
		\end{equation}
		Let $A:=\begin{pmatrix}
			1&0&0&0\\
			0&1&0&0\\
			0&0&\frac{1}{\sqrt{2}}&-\frac{1}{\sqrt{2}}\\
			0&0&\frac{1}{\sqrt{2}}&\frac{1}{\sqrt{2}}
		\end{pmatrix}$ and $\pi:=\Xi^{-1}\circ \tilde{\pi}\circ A: S^3\to U^*S^2$ so that 
		\begin{equation}\label{eq:pi_inv}
			\pi^{-1}(\Gamma\cup U^*_nS^2)\subset S^3\cap \mathrm{span}\{1,\mathbf{i},\mathbf{j}\}.
		\end{equation} 
		As pointed out in \cite[Chapter 1.1]{massot2014topological}, the map $\Xi^{-1}\circ\tilde\pi$ pulls back the contact structure on $U^*S^2$ to the standard one on $S^3$. The same holds for $\pi$ since  $A$ is symplectic.
	\end{proof}
	
	\begin{remark}
		Note that for $U^*S^2\# U^*S^2$ in \eqref{eq:cont_conn} two copies of $U^*S^2$ are glued near the tangent vector $(-1,0,0)$ at the south pole, see \eqref{eq:lagrange}. This tangent vector corresponds to $\begin{pmatrix}
			-1&0&0\\
			0&0&-1\\
			0&-1 &0
		\end{pmatrix}$, which is the rotation by angle $\pi$ with respect to the rotation axis spanned by $(0,-1,1)$. Thus, it lifts via $\pi$ to $\pm\cos(\frac{\pi}{2})\mp\sin (\frac{\pi}{2})\mathbf{k}=\mp\mathbf{k}$, the south and north poles of $S^3\subset\text{Im}(\mathbb{H})$. 
	\end{remark}

	Let $(x,y,z)$ be the coordinates of $\mathbb{R}^3$, and let $\eta$ be the angular coordinate on $S^1$. We consider the 1-form 
	\[
	\lambda= z\text{d}\eta+x\text{d}y-y\text{d}x 
	\]
	on $ S^1\times \mathbb{R}^3$ such that $\D\lambda$ is a symplectic form. The involution 
	\begin{equation}\label{eq:sym}
		\sigma: S^1\times \mathbb{R}^3\to S^1\times \mathbb{R}^3, \qquad (\eta,x,y,z)\mapsto (-\eta,-x,-y,-z)
	\end{equation}
	preserves $\lambda$ and therefore is symplectic. The hypersurface $S^1\times S^2$ in  $S^1\times \mathbb{R}^3$, where $S^2$ denotes the unit sphere, is of contact type with $\alpha:=\lambda|_{S^1\times S^2}$, and $\sigma$ restricted to $S^1\times S^2$ is an involution preserving $\alpha$.
	
	Let $\alpha^\mathfrak{B}$ be the contact form on $\overline{\Sigma}_d^\mathfrak{B}$ in Proposition \ref{prop:contact_birkhoff}, and let $\sigma^{\mathfrak{B}}$ be the nontrivial deck transformation of  $\pi^{\mathfrak{B},\mathfrak{M}}:\overline{\Sigma}_d^{\mathfrak{B}} \to \overline{\Sigma}_d^{\mathfrak{M}}$ in \eqref{eq:cover_moser}. Let us denote $\xi^\mathfrak{B}:=\ker\alpha^\mathfrak{B}$. For the standard contact form $\alpha$ on $S^1\times S^2$ introduced above, we denote $\xi:=\ker\alpha$. 
	\begin{proposition}
		For every $d\in(H(L_1),H(L_1)+\epsilon)$, there exists a contactomorphism 
		\[
		F:(\overline{\Sigma}_d^\mathfrak{B},\xi^\mathfrak{B}=\ker\alpha^\mathfrak{B}) \to (S^1\times S^2,\xi=\ker\alpha)
		\]
		such that $F(\Lambda_d^{\mathfrak{B},E})$, $F(\Gamma_d^{\mathfrak{B},E})$ are meridians in  $\{0\}\times S^2$ and $F(\Lambda_d^{\mathfrak{B},M})$, $F(\Gamma_d^{\mathfrak{B},M})$ are meridians in $\{\pi\}\times S^2$. Furthermore, it satisfies $F\circ \sigma^{\mathfrak{B}}= \sigma\circ F$.	
		\label{cansubst}
	\end{proposition}
	\begin{proof}
		In the hyperspherical coordinates on $S^3\subset\mathbb{R}^4$, namely
		\[
		(x_1,x_2,x_3,x_4)=\big(\sin(\eta)\sin(\theta)\cos(\varphi),\sin(\eta)\sin(\theta)\sin(\varphi),\sin(\eta)\cos(\theta),\cos(\eta)\big),
		\]
		the standard contact form on $S^3$ is written as
		\begin{align*}
			x_1\text{d}x_2-x_2\text{d}x_1+x_3\text{d}x_4-x_4\text{d}x_3=\cos(\theta)\text{d}\eta-\sin(\eta)\cos(\eta)\sin(\theta)\text{d}\theta+\sin^2(\eta)\sin^2(\theta)\text{d}\varphi.
		\end{align*}
		We now consider two copies of $S^3$ and perform a connected-sum operation twice. More precisely,  we take two copies of $S^3$, each with small open balls removed around the north and south poles; that is, we consider
		\[
		\mathring S^3:= S^3\setminus (\{0\leq \eta \leq \delta\}\cup \{\pi-\delta\leq \eta\leq\pi\})=\{\delta<\eta<\pi-\delta\}\times S^2
		\] 
		for small $\delta>0$. We then smoothly connect the boundaries near the north poles of each $\mathring S^3$ using $[0,1]\times S^2$, and do the same for the south poles.  We denote the resulting smooth manifold by $S^3\#_2 S^3$. The covering map $S^3\to U^*S^2$ in Lemma \ref{lem:covering} induces a double covering map 
		\[
		\pi: S^3\#_2 S^3\longrightarrow U^*S^2\# U^*S^2.
		\]
		Moreover, we have the following commutative diagram:
		\[
		\begin{tikzcd}		
			\big(\overline{\Sigma}_d^\mathfrak{B},\xi^\mathfrak{B}\big) \arrow{r}{G}\arrow{d}{\pi^{\mathfrak{B},\mathfrak{M}}} &\big(S^3\#_2 S^3,\pi^* {\xi}^\#\big) \arrow{d}{\pi}\\
			\big(\overline{\Sigma}_d^\mathfrak{M}, \xi^\mathfrak{M}\big) \arrow{r}{g} & \big(U^*S^2\# U^*S^2, {\xi}^\#\big) 
		\end{tikzcd}   
		\]
		where $\xi^\mathfrak{M}$ is the contact structure given in Theorem \ref{thm:moser2}. By the definition of $\xi^\mathfrak{B}$, we have $(\pi^{\mathfrak{B},\mathfrak{M}})^*\xi^\mathfrak{M}=\xi^\mathfrak{B}$. The contactomorphism $g$ is described in \eqref{eq:cont_conn}, where $\xi^\#$ is the contact structure given by performing the connected sum operation to two copies of the standard contact $U^*S^2$. The map $G$ is the lift of $g$. Thus, it maps $\xi^\mathfrak{B}$ to $\pi^* {\xi}^\#$ and commutes with the deck transformations of $\pi^{\mathfrak{B},\mathfrak{M}}$ and $\pi$. 
		
		Now we construct a diffeomorphism 
		\[
		h:S^3\#_2 S^3\longrightarrow S^1\times S^2
		\]
		in the following way. Using the subscripts $1,\,2$ to distinguish the first and second $S^3$ of $S^3\#_2 S^3$, we map 
		\begin{equation*}
			\begin{split}
				\mathring S_1^3=\{\delta< \eta_1< \pi-\delta\}\times S^2&\longrightarrow S^1\times S^2\\
				(\eta_1,x)&\longmapsto (\eta_1+\tfrac{\pi}{2},x),
			\end{split}
		\end{equation*}
		and
		\begin{equation*}
			\begin{split}
				\mathring S_2^3=\{\delta< \eta_2< \pi-\delta\}\times S^2&\longrightarrow S^1\times S^2  \\
				(\eta_2,x)&\longmapsto (\tfrac{5\pi}{2}-\eta_2,x).
			\end{split}
		\end{equation*}
		We then extend this map to the entire $S^3_1\#_2S^3_2$ in such a way that two copies of $[0,1]\times S^2$, which connect $\mathring S^3_1$ and $\mathring S^3_2$, are mapped to the remaining part  of $S^1\times S^2$, and the extended map $h$ is a diffeomorphism equivariant with respect to $\sigma$  in \eqref{eq:sym} and the nontrivial deck transformation of $\pi$. Remembering that $\pi^*\xi^\#$ is the standard contact structure on $S^3$ away from the connected sum region by Lemma \ref{lem:covering}, we verify that  $h_*\pi^*\xi^\#$ coincides with $\xi$ on  $\{0\}\times S^2=h(\{\eta_2=\frac{\pi}{2}\})$ and $\{\pi\}\times S^2=h(\{\eta_1=\frac{\pi}{2}\})$:
		\begin{align*}
			\cos(\theta)\text{d}\left(\eta-\frac{\pi}{2}\right)-\sin\left(\frac{\pi}{2}\right)\cos\left(\frac{\pi}{2}\right)\sin(\theta) \text{d}\theta+\sin^2\left(\frac{\pi}{2}\right)\sin^2(\theta)\text{d}\varphi = \cos(\theta)\text{d}\eta+\sin^2(\theta)\text{d}\varphi.
		\end{align*}

		The contact structure $h_*\pi^*\xi^\#$ on $S^1\times S^2$ is tight. We now construct a contactomorphism 
		\begin{equation}\label{eq:cont_j}
			j: (S^1\times S^2, h_*\pi^*\xi^\#)\to (S^1\times S^2, \xi)	
		\end{equation}
		which is the identity on $\{0,\pi\}\times S^2$.	By Giroux's theorem, see \cite[Theorem~2.5.22 and Theorem~2.5.23]{geiges2008a}, we have a diffeomorphism  
		\[
		j_0:\left([-\epsilon,\epsilon]\times S^2\right)\cup \left([\pi-\epsilon,\pi+\epsilon]\times S^2\right) \to \left([-\epsilon,\epsilon]\times S^2\right)\cup \left([\pi-\epsilon,\pi+\epsilon]\times S^2\right)
		\] 
		such that $j_0^*\xi=h_*\pi^*\xi^\#$, it is the identity map on $\{0,\pi\}\times S^2$, and $j_0\circ\sigma=\sigma\circ j_0$. Due to  \cite[Theorem~4.9.4]{geiges2008a} in combination with \cite[Remark~4.9.3]{geiges2008a}, we can extend $j_0$ to a contactomorphism on $[-\epsilon,\pi+\epsilon]\times S^2$. Using the involution $\sigma$, we further extend it over $[\epsilon,\pi-\epsilon]\times S^2$. This yields the contactomorphism claimed in \eqref{eq:cont_j}, which is clearly equivariant with respect to $\sigma$.  
		
		The map $F:=j\circ h\circ G$ is the equivariant contactomorphism claimed in the statement of the proposition. By  Lemma \ref{lem:covering} and the construction of $F$,  $F(\Lambda_d^{\mathfrak{B},E})\cup F(\Gamma_d^{\mathfrak{B},E}) \subset \{0\}\times S^2$ and $F(\Lambda_d^{\mathfrak{B},M})\cup F(\Gamma_d^{\mathfrak{B},M}) \subset \{\pi\}\times S^2$. Note that all Legendrian knots in $\{0,\pi\}\times S^2$ are necessarily meridians since $\alpha|_{\{0,\pi\}\times S^2}=\sin^2(\theta)d\varphi$, and this completes the proof.
	\end{proof}

	\subsection{Indices of Reeb chords on $S^1\times S^2$}\label{ReebChords}
	In this section, we will study Reeb chords of $(S^1\times S^2, \alpha)$ whose endpoints lie on a meridian of $\{0\} \times S^2$. We consider a meridian
	\[
	\Lambda_{0}:=\left\{(\eta=0, \varphi=0, \theta)\ |\ \theta\in [0,2\pi] \right\}\subset \{0\}\times S^2,
	\]
	which is a Legendrian knot in $S^1\times S^2$. 
	Remember that the standard contact form on $S^1\times S^2$ is given by
	\[
	\lambda= z\text{d}\eta+\big(x\text{d}y-y\text{d}x), 	 
	\]
	where $(x,y,z)$ are the coordinates of $\mathbb{R}^3$ and we view $S^2$ as an sphere embedded in $\mathbb{R}^3$. Using the spherical coordinates
	\[
	(x,y,z)=(\sin\theta \cos\varphi, \sin\theta\sin\varphi,\cos\theta),
	\]
	we write 
	\[
	\lambda= \cos(\theta)\text{d}\eta +\sin^2(\theta)\text{d}\varphi,\quad 	\text{d}\lambda= -\sin(\theta)\text{d}\theta\wedge \text{d}\eta +2\sin(\theta)\cos(\theta)\text{d}\theta\wedge\text{d}\varphi.
	\]
	The Reeb vector field is given by
	\[
	R(\eta, \varphi, \theta)= \frac{2\cos(\theta)}{1+\cos^2(\theta)}\del_\eta + \frac{1}{1+\cos^2(\theta)}\del_\varphi
	\]
	and the corresponding Reeb flow is
	\[	\Phi_R^t( \eta_0, \theta_0,\varphi_0) = 
	\Bigr(
	\frac{2\cos(\theta_0)}{1+\cos^2(\theta_0)}t+\eta_0 ,\, 
	\theta_0,\, \frac{1}{1+\cos^2(\theta_0)}t+\varphi_0 \Bigr)
	\]
	We are interested in Reeb chords that begin and end in $\Lambda_{0}$ and that have the trivial homotopy class in $\pi_1(S^1\times S^2,\Lambda_0)$. They are 
	\begin{align*}
		&\gamma^k_1(t):= \left(\eta=0,\varphi=t, \theta=\frac{\pi}{2}\right)\ \text{ for }\ t\in [0,k\pi],\\[1ex]
		&\gamma^k_2(t):= \left(\eta=0,\varphi=t+\pi, \theta=\frac{\pi}{2}\right)\ \text{ for }\ t\in [0,k\pi].
	\end{align*}
	Note that $\pi_1(S^1\times S^2,\Lambda_0)=\pi_1(S^1\times\mathbb{R}^3,L_0)$, where $L_0$ is the $xz$-plane in $\mathbb{R}^3$. 
	
	\medskip
	
	To compute the Robbin-Salamon indices for $\gamma^k_1$ and $\gamma^k_2$, we now think of them as chords in $S^1\times\mathbb{R}^3$ with endpoints in $L_0$. Let $r$ denote the radial coordinate on $\mathbb{R}^3$. As usual, we extend the Reeb flow in the radial direction trivially, i.e.~
	\[
	\Phi_R^t ( \eta_0,r_0, \theta_0,\varphi_0 ) = 
	\Bigg(
	\frac{2\cos(\theta_0)}{1+\cos^2(\theta_0)}t+\eta_0 ,r_0, \theta_0,
	\frac{1}{1+\cos^2(\theta_0)}t+\varphi_0
	\Bigg).
	\]
	In the frame $\{\del_\eta, \del_r,\del_\theta, \del_\varphi\}$, the linearized flow of $\Phi_R^t$ along $\gamma_1^k$ is as follows:
	\[
	\text{D}\Phi^t_R(\eta, r, \theta, \varphi) = \begin{pmatrix}
		1&0& \frac{2\sin(\theta)\left(\cos^2(\theta)-1\right) }{\left(1+\cos^2(\theta)\right)^2}t &0\\
		0&1&0&0\\
		0&0&1& 0\\
		0&0&\frac{2\cos(\theta)\sin(\theta)}{\left(1+\cos^2(\theta)\right)^2}t &1
	\end{pmatrix},\quad 
	\text{D}\Phi^t_R\left(\gamma_1^k(0)\right)  = \begin{pmatrix}
		1&0&-2t&0\\
		0&1&0&0\\
		0&0&1&0\\
		0&0&0&1
	\end{pmatrix}.
	\]
	Since the frame $\{\del_\eta, \del_r,\del_\theta, \del_\varphi\}$ is not globally well-defined, we now switch to the Cartesian frame $\{\del_x,\del_
	z,\del_y,\del_\eta\}$, with the basis ordered so that the tangent space of $L_0$ at $\gamma_1^k(0)$ becomes $\mathbb{R}^2\times \{(0,0)\}$. Recall that $\D\lambda=2\D x\wedge \D y + \D z\wedge\D \eta$. In this Cartesian frame, $\text{D}\Phi^t_R\left(\gamma_1^k(0)\right)$ is written as 
	\[
	\Psi(t):=\begin{pmatrix}
		\cos(t)&0&-\sin(t)&0\\
		0&1&0&0\\
		\sin(t)&0&\cos(t)&0\\
		0&2t&0&1
	\end{pmatrix}.
	\]
	We write $\Psi(t):= \begin{pmatrix} A(t) & B(t) \\ C(t) & D(t)		\end{pmatrix}$, then the Robbin-Salamon index of $\gamma_1^k$ is computed as
	\[
	\mu_{RS}(\gamma_1^k)=\frac{1}{2}\mathrm{sign\,} \Gamma(0)
	+\sum_{0<t<k\pi}\mathrm{sign\,}\Gamma(t) + \frac{1}{2}\mathrm{sign\,} \Gamma(k\pi),
	\]
	where $\mathrm{sign}$ denotes the signature and 
	\[
	\Gamma(t):\ker C(t)\to\mathbb{R};\quad \Gamma(t)v=\langle A(t)v, P
	\dot  C(t)v\rangle.
	\]
	Here $P=\begin{pmatrix} 2& 0\\ 0& 1	\end{pmatrix}$ appears due to the coefficient 2 in $\D\lambda$. We refer to \cite[Remark 2.5]{robbin1993a} for details. Then, a straightforward computation shows that 
	\begin{equation}
		\mu_{RS}(\gamma_1^k)=\mu_{RS}(\gamma_2^k)= k+\frac{1}{2}, \qquad\forall k\in \mathbb{N},
		\label{index1}
	\end{equation}
	where $\mu_{RS}(\gamma_1^k)=\mu_{RS}(\gamma_2^k)$ holds due to symmetry.
	
	\medskip

	\section{Consecutive Collision Orbits and Floer homology}
	\label{orbits}
	In this section, we will compute the Lagrangian Rabinowitz Floer homology of $(S^1\times S^2,\alpha)$ with respect to Legendrian knots given by meridians in $\{0\}\times S^2$. Using this computation, we then show the existence of consecutive collision orbits in the R3BP for energy levels slightly above the first critical energy value.
	
	\subsection{Equivariant Lagrangian RFH}\label{sec:equiv}
	
	Let us endow $S^1\times\mathbb{R}^3$ with the symplectic form $\D\lambda$, where $\lambda=z \D\eta+x\D y- y \D x$, and let us choose a Lagrangian submanifold $L_0:=\{0\}\times\mathbb{R}\times \{0\}\times\mathbb{R}$. We consider a contact form $e^f\alpha$ on $S^1\times S^2$, which supports the same contact structure $\xi$ as $\alpha$, and choose an embedding $\iota_f:S^1\times S^2\hookrightarrow S^1\times\mathbb{R}^3$ so that $\iota_f^*\lambda=e^f\alpha$. We denote 
	\[
	\Sigma:=\iota_f(S^1\times S^2),\qquad \Lambda^\Sigma_0:=\Sigma\cap L_0.
	\]
	The latter is a Legendrian knot in $(\Sigma,\lambda|_{\Sigma})$. From now on, we assume that every Reeb chord on $\Sigma$ with endpoints on $\Lambda_0^\Sigma$ is nondegenerate. 
	
	To introduce the Rabinowitz action functional associated to $(\Sigma,\Lambda_0^\Sigma)$, let $H:S^1\times\mathbb{R}^3\to \mathbb{R}$ be a smooth function such that $H^{-1}(0)=\Sigma$ is a regular level set, the Hamiltonian vector field $X_H$ defined by $\iota_{X_H} d\lambda= -dH$ coincides with the Reeb vector field on $(\Sigma,\lambda|_\Sigma)$, and $H$ is constant outside a bounded region. Let $P\left(S^1\times \mathbb{R}^3,L_{0};[\text{pt.}]\right)$ be the space of smooth paths $\gamma:([0,1],\{0,1\})\to(S^1\times \mathbb{R}^3,L_0)$ having the trivial homotopy class in $\pi_1(S^1\times\mathbb{R}^3,L_0)$. The Rabinowitz action functional is defined by
	\[
	\mathcal{A}^H:\; P\left(S^1\times \mathbb{R}^3,L_{0};[\text{pt.}]\right)\times\mathbb{R}\longrightarrow\mathbb{R},\qquad\mathcal{A}^H(\gamma,\tau):=\int\limits_0^1 \gamma^*\lambda - \tau \int\limits_0^1 H(\gamma(t))\D t.
	\]
	A pair $(\gamma,\tau)$ is a critical point of $\mathcal{A}^H$ if and only if 
	\[
	\del_t \gamma(t)=\tau X_H(\gamma(t)),\qquad \gamma(t)\in \Sigma,\qquad \forall t\in [0,1].
	\]
	The former implies that $\bar\gamma(t):=\gamma(t/\tau)$ for $t\in [0,\tau]$ is a (generalized) Reeb chord with  both endpoints on $\Lambda_0^\Sigma$. By generalized, we mean that $\tau$ can be nonpositive. If $\tau=0$, then $\gamma$ is a constant chord on $\Lambda_0^\Sigma$. If $\tau<0$, then $\bar\gamma$ is an orbit of the negative Reeb vector field. The connected components of $\mathrm{Crit\,}\mathcal{A}^H$ that consist of constant chords is diffeomorphic to $\Lambda_0^\Sigma$. We choose a Morse function $h$ on this space, which is diffeomorphic to a circle, with maximum points $(\gamma_{\mathrm{const},\ell}^+ ,0)$ and minimum points $(\gamma_{\mathrm{const},\ell}^- ,0)$, where $\ell\in\{1,\dots,\ell_0\}$ for some $\ell_0\geq1$. The Lagrangian Rabinowitz Floer chain modules are defined by
	\[
	RFC_i\left(S^1\times\mathbb{R}^3, \Sigma, L_{0};[\text{pt.}]\right)=\begin{cases} 
		\displaystyle\bigoplus_{(\gamma,\tau)}\mathbb{Z}_2\langle(\gamma,\tau)\rangle \qquad & i\in\mathbb{Z}\setminus\{0,1\}\\[2ex]
		\displaystyle\bigoplus_{(\gamma,\tau)}\mathbb{Z}_2\langle(\gamma,\tau)\rangle \oplus \bigoplus_{1\leq \ell\leq\ell_0} \mathbb{Z}_2\langle(\gamma_{\mathrm{const},\ell}^+ ,0)\rangle & i=1\\[2ex]
		\displaystyle\bigoplus_{(\gamma,\tau)}\mathbb{Z}_2\langle(\gamma,\tau)\rangle \oplus \bigoplus_{1\leq \ell\leq\ell_0}\mathbb{Z}_2\langle(\gamma_{\mathrm{const},\ell}^- ,0)\rangle & i=0
	\end{cases}
	\] 
	where each direct sum $\displaystyle\bigoplus_{(\gamma,\tau)}$ ranges over all critical points $(\gamma,\tau)$ of $\mathcal{A}^H$ such that $\tau\neq 0$ and $\mu_{RS}(\bar\gamma)+\frac{1}{2}=i$. As usual, the boundary operator 
	\[
	\partial_i: RFC_i\left(S^1\times\mathbb{R}^3, \Sigma, L_{0};[\text{pt.}]\right)\to  RFC_{i-1}\left(S^1\times\mathbb{R}^3, \Sigma, L_{0};[\text{pt.}]\right)
	\] 
	is defined by counting $L^2$-gradient flow lines of $\mathcal{A}^H$ (together with gradient flow lines of the Morse function $h$ in degree 0 and 1).
	We write $RFH_*\left(S^1\times\mathbb{R}^3, \Sigma, L_{0};[\text{pt.}]\right)$ for the homology of this chain complex. We refer to \cite{merry2014a} for details on the construction of the Lagrangian Rabinowitz Floer homology. 
	\medskip
	
	Since the hypersurface $\Sigma$ can be displaced from the Lagrangian $L_{0}$ by a compactly supported Hamiltonian in $S^1\times \mathbb{R}^3$, due to \cite[Chapter 2.4]{merry2014a}, we have 
	\begin{equation}\label{eq:vanishing}
		RFH_i\left(S^1\times\mathbb{R}^3, \Sigma, L_{0};[\text{pt.}]\right)=0 \qquad \forall i\in\mathbb{Z}.
	\end{equation}
	Since there are constant generators $(\gamma_\mathrm{const,\ell}^\pm ,0)$, the above vanishing result implies the existence of a nonconstant Reeb chord on $(\Sigma,\Lambda_0^\Sigma)$. To go beyond the mere existence result, we incorporate the symplectic involution $\sigma$ defined in \eqref{eq:sym} into our analysis.
	
	From now on, we assume that $\Sigma$ is invariant under $\sigma$. We may assume that the Morse function $h$ is invariant under $\sigma$. Since $\sigma^*\lambda=\lambda$ and $\sigma(L_0)=L_0$, $\sigma$ acts freely on the chain complex $RFC_i\left(S^1\times\mathbb{R}^3, \Sigma, L_{0};[\text{pt.}]\right)$. The boundary operators $\partial_i$ can be made equivariant with respect to the action induced by $\sigma$ by defining them with a $\sigma$-equivariant almost complex structure.  
	Therefore, we can define the $\mathbb{Z}_2$-equivariant complex and denote its homology by
	\[
	RFH_i^{\mathbb{Z}_2}\left(S^1\times\mathbb{R}^3, \Sigma, L_{0};[\text{pt.}]\right),
	\] 
	where the superscript $\mathbb{Z}_2$ refers to the $\mathbb{Z}_2$-action generated by $\sigma$, cf.~\cite[Theorem~4.5.6]{ruck24b}. This equivariant homology can be computed by means of the Tate homology, see \cite{ruck2024a}. Below we provide a direct  computation.

	\begin{proposition}\label{prop:equiv_rfh}
		Let $\Sigma$ be as above. Then,
		\[
		RFH_i^{\mathbb{Z}_2}\left(S^1\times\mathbb{R}^3, \Sigma, L_{0};[\text{pt.}]\right)\cong \mathbb{Z}_2\qquad  \forall i\in \mathbb{Z}.
		\]
	\end{proposition}
	\begin{proof}
		A standard continuation argument shows that 
		\[
		RFH_i^{\mathbb{Z}_2}\left(S^1\times\mathbb{R}^3, \Sigma, L_{0};[\text{pt.}]\right)\cong RFH_i^{\mathbb{Z}_2}\left(S^1\times\mathbb{R}^3, S^1\times S^2, L_{0};[\text{pt.}]\right).
		\]
		It suffices to show
		\begin{equation}\label{eq:rfh_s1xs2}
			RFH_i^{\mathbb{Z}_2}\left(S^1\times\mathbb{R}^3, S^1\times S^2, L_{0};[\text{pt.}]\right)\cong \mathbb{Z}_2\qquad \forall i\in\mathbb{Z}.
		\end{equation}
		While this computation can be carried out using Tate homology, it can also be seen directly as follows. We identify the space of constant Reeb chords on $(S^1\times S^2,\alpha)$ with  $\Lambda_0:=(S^1\times S^2)\cap L_0$. We choose a Morse function $h$ on $\Lambda_0$ that is invariant under $\sigma|_{\Lambda_0}$ and has exactly two maximum points $\gamma_\mathrm{const,1}^+,\,\gamma_\mathrm{const,2}^+$ and two minimum points $\gamma_\mathrm{const,1}^-,\,\gamma_\mathrm{const,2}^-$. By the index computation in \eqref{index1}, we have
		\begin{equation}\label{eq:rfc_s1xs2}
			RFC_i\left(S^1\times\mathbb{R}^3, S^1\times S^2, L_{0};[\text{pt.}]\right)=\begin{cases} 
				\mathbb{Z}_2\big\langle(\gamma^i_1, i\pi),(\gamma^i_2, i\pi)\big\rangle \qquad & i\in\mathbb{Z}\setminus\{0,1\}\\[2ex]
				\mathbb{Z}_2\big\langle(\gamma_\mathrm{const,1}^+ ,0),(\gamma_\mathrm{const,2}^+ ,0)\big\rangle & i=1\\[2ex]
				\mathbb{Z}_2\big\langle(\gamma_\mathrm{const,1}^- ,0),(\gamma_\mathrm{const,2}^- ,0)\big\rangle & i=0\,.
			\end{cases}	
		\end{equation}
		Note that the $\mathbb{Z}_2$-action induced by $\sigma$ maps one generator to the other one in each degree. Moreover, the boundary operators $\partial_i$ are equivariant with respect to the $\mathbb{Z}_2$-action. 
		By \eqref{eq:rfc_s1xs2}, for a degree reason, $\partial_1$ is just the Morse boundary operator for $h$ on $\Lambda_0\cong S^1$. Therefore, we know that
		\[
		\partial_1(\gamma_\mathrm{const,1}^+,0)=\partial_1(\gamma_\mathrm{const,2}^+,0)= (\gamma_\mathrm{const,1}^-,0)+(\gamma_\mathrm{const,2}^-,0)\,.
		\]
		The facts that $\partial_i\circ \partial_{i+1}=0$ for all $i\in\mathbb{Z}$ and that the chain complex in \eqref{eq:rfc_s1xs2} is acyclic by \eqref{eq:vanishing} imply that 
		\[
		\partial_k(\gamma_1^k,k\pi)=\partial_k(\gamma_2^k,k\pi)=(\gamma_1^{k-1}, (k-1)\pi)+(\gamma_2^{k-1},(k-1)\pi) \qquad \forall k\in\mathbb{Z}.
		\]
		Therefore, the equivariant complex $RFC_i^{\mathbb{Z}_2}\left(S^1\times\mathbb{R}^3, S^1\times S^2, L_{0};[\text{pt.}]\right)$ has rank one in each degree, and the boundary operators in this equivariant complex vanish. Note that we use $\mathbb{Z}_2$-coefficients. This proves \eqref{eq:rfh_s1xs2}.
	\end{proof}

	\subsection{Perturbed LRFH and symmetric consecutive collision orbits}
	\label{sec:perturb_RFH}
	Proposition \ref{prop:equiv_rfh} implies the existence of infinitely many Reeb chords on $\Sigma$ that begin and end in $\Lambda_0^\Sigma$. However, if there is a Reeb chord that is also periodic, it produces infinitely many Reeb chords through iterations. This issue will be partially addressed in the context of the R3BP in Proposition \ref{noP}. Before that, we introduce a perturbed version of Rabinowitz Floer homology that can be used to detect Reeb chords between two Hamiltonian isotopic Legendrians.

	We continue to assume that $\Sigma\subset S^1\times \mathbb{R}^3$ is invariant under the symplectic involution $\sigma$. Let us denote by
	\[
	L_{\varphi_0}:=\{0\}\times \{(r\sin\theta \cos\varphi_0, r\sin\theta\sin\varphi_0,r\cos\theta) \mid 0\leq r,\,0\leq\theta<2\pi\}\subset S^1\times\mathbb{R}^3
	\]
	the plane obtained by rotating the $xz$-plane by an angle $\varphi_0\in [0,\pi)$. We seek  Reeb chords that start on $\Lambda_0^\Sigma$ and end on $\Lambda_{\varphi_0}^\Sigma:=\Sigma \cap L_{\varphi_0}$.
	
	We use a perturbed version of the Rabinowitz action functional, originally introduced in  \cite{albers2010leaf} to detect leaf-wise intersections. The Lagrangian analogue, which we employ here, is due to \cite{merry2014a}. We consider the Hamiltonian 
	\[
	G: S^1\times \mathbb{R}^3\to \mathbb{R},\qquad  (\eta, x,y,z)\mapsto \varphi_0 \left(x^2+ y^2\right).
	\]
	The time one-map $\phi_G^1$ of its Hamiltonian flow maps $L_{\varphi_0}$ to $L_0$. By multiplying with a suitable cutoff function, we may assume that $G$ is compactly supported and 
	\begin{equation}\label{eq:time_one_map_G}
		\phi_G^{-1}(L_0)\cap \Sigma=\Lambda^\Sigma_{\varphi_0}	
	\end{equation} 
	We further modify $G$ to  depend on time in such a way that it is 1-periodic in time with support in the time interval $(\frac{1}{2},1)$, while ensuring that the time-one map of the Hamiltonian flow remains unchanged, see \cite[Lemma 2.3]{albers2010leaf}. Abusing notation, we continue to denote the modified function by $G$. We now define the perturbed Rabinowitz action functional
	\[
	\begin{aligned}
		&\mathcal{A}_G^H: P\left(S^1\times \mathbb{R}^3,L_{0};[\text{pt.}]\right)\times\mathbb{R}\longrightarrow \mathbb{R},\\
		&\mathcal{A}_G^H(x,\tau):=\int\limits_0^1 x^*\lambda  - \tau \int\limits_0^1\beta(t)H(x(t))\D t - \int\limits_0^1 G(x(t),t)\D t,
	\end{aligned}
	\]
	where $H:S^1\times\mathbb{R}^3\to\mathbb{R}$ is the smooth function associated to the hypersurface $\Sigma$ used in the previous section, and  $\beta:[0,1]\to\mathbb{R}$ is a smooth function supported in $\left(0,\frac{1}{2}\right)$ with total integral equal to $1$. 	 A critical point of this functional is  a pair $(x,\tau)$ with
	\[
	\begin{aligned}
		\del_t x(t)=\tau\beta(t)X_H(x(t))&+X_G(x(t),t),\\
		\int\limits_0^1 \beta(t)H(x(t))\D t&=0.
	\end{aligned}
	\]
	Recall that $X_H$ and $X_G$ have disjoint support in time, namely $(0,\frac{1}{2})$ and $(\frac{1}{2},1)$, respectively. One can readily see that a critical point $(x,\tau)$ of $\mathcal{A}_G^H$ corresponds to a relative leaf-wise intersection point of $\phi_G^1$, meaning that 
	\[
	x(0)\in\Lambda^\Sigma_0=\Sigma\cap L_0,\qquad \phi_G^1(\phi_R^\tau(x(0)))\in L_0.
	\] 
	The latter property together with \eqref{eq:time_one_map_G} implies $\phi_R^\tau(x(0))\in \Lambda^\Sigma_{\varphi_0}$. Hence, a critical point $(x,\tau)$  gives rise to a Reeb chord $\gamma(t):=\phi_R^t(x(0))$ on $\Sigma$ that starts at $\gamma(0)\in\Lambda^\Sigma_0$ and ends at $\gamma(\tau)\in \Lambda^\Sigma_{\varphi_0}$. 
	
	\medskip
	We now consider the chain complex $RFC_i\left(S^1\times\mathbb{R}^3,\Sigma, L_{0}, G, [\text{pt.}]\right)$, which are $\mathbb{Z}_2$-modules generated by the critical points of $\mathcal{A}_G^H$ and graded by the Robbin-Salamon index. The boundary operator is defined by counting $L^2$-gradient flow lines of $\mathcal{A}_G^H$.
	We may assume that $G$ is invariant under the symplectic involution $\sigma$ in \eqref{eq:sym} since the original form of $G$ prior to multiplication by a cutoff function is invariant.  We can therefore  define the $\mathbb{Z}_2$-equivariant chain complex $RFC^{\mathbb{Z}_2}_i\left(S^1\times\mathbb{R}^3,\Sigma, L_{0}, G, [\text{pt.}]\right)$ and its homology $RFH^{\mathbb{Z}_2}\left(S^1\times\mathbb{R}^3, \Sigma , L_{0}, G, [\text{pt.}]\right)$. A $\sigma$-invariant homotopy between $G$ and the zero function induces a continuation isomorphism 
	\begin{equation}\label{RFH3}
		RFH_i^{\mathbb{Z}_2}\left(S^1\times\mathbb{R}^3, \Sigma , L_{0}, G, [\text{pt.}]\right)\cong RFH_i^{\mathbb{Z}_2}\left(S^1\times\mathbb{R}^3, \Sigma, L_{0}, [\text{pt.}]\right)\cong \mathbb{Z}_2
	\end{equation}
	for every $i\in\mathbb{Z}$, where the last isomorphism was proved in Proposition \ref{prop:equiv_rfh}.

	\medskip

	We now apply our general discussion to establish the existence of symmetric consecutive collision orbits in the R3BP for energy levels slightly above the first critical energy value. 
	
	\begin{corollary}\label{cor2}
		Let $\epsilon>0$ be as in Theorem \ref{thm:moser2}. Then, for every  $d\in(H(L_1),H(L_1)+\epsilon)$, the unregularized energy hypersurface $\left(\Sigma_d,X_{H}\right)$ admits either a symmetric periodic collision orbit or infinitely many consecutive collision orbits with the Earth. The same statement holds for the Moon.
	\end{corollary}
	\begin{proof}
		Recall from Proposition~\ref{cansubst} that we have a contactomorphism
		\[
		F:(\overline{\Sigma}_d^\mathfrak{B},\xi^\mathfrak{B}=\ker\alpha^\mathfrak{B}) \longrightarrow (S^1\times S^2,\xi=\ker\alpha).
		\]
		Thus, $(F^{-1})^*\alpha^{\mathfrak B}=e^{f_{\mathfrak{B}}}\alpha$ for some $f_{\mathfrak{B}}\in C^\infty(S^1\times S^2)$. Thus, there is a hypersurface in $S^1\times\mathbb{R}^3$ corresponding the contact form $e^{f_{\mathfrak{B}}}\alpha$. Abusing notation, we denote this hypersurface by $\overline{\Sigma}_d^\mathfrak{B}$. Due to $F\circ \sigma^{\mathfrak{B}}= \sigma\circ F$ proved in Proposition~\ref{cansubst}, we know $\overline{\Sigma}_d^\mathfrak{B}\subset S^1\times \mathbb{R}^3$ is invariant under $\sigma$. Recall also that  $F(\Lambda_d^{\mathfrak{B},E})$ and $F(\Gamma_d^{\mathfrak{B},E})$ are meridians in  $\{0\}\times S^2$. Therefore 		
		\begin{equation}\label{eq:leg}
			\overline{\Sigma}_d^\mathfrak{B}\cap L_0 \quad \text{and}\quad  \overline{\Sigma}_d^\mathfrak{B}\cap L_{\varphi_0}
		\end{equation}
		for some $\varphi_0\in[0,\pi)$  correspond to $\Lambda_d^{\mathfrak{B},E}$ and $\Gamma_d^{\mathfrak{B},E}$, respectively. Let $G:S^1\times \mathbb{R}^3\to\mathbb{R}$ be as above so  that it has the property $\phi_G^{-1}(L_0)\cap \overline{\Sigma}_d^\mathfrak{B}=\overline{\Sigma}_d^\mathfrak{B}\cap L_{\varphi_0}$, cf.~\eqref{eq:time_one_map_G}. Assume that $(\overline{\Sigma}_d^\mathfrak{B},\Lambda_d^{\mathfrak{B},E},G)$ is nondegenerate so that the associated homologies are well-defined. Then, by \eqref{RFH3}, we have
		\[
		RFH_i^{\mathbb{Z}_2}\big(S^1\times\mathbb{R}^3, \overline{\Sigma}_d^\mathfrak{B} , L_{0}, G, [\text{pt.}]\big)\cong \mathbb{Z}_2.
		\]
		Thus, the associated action functional $\mathcal{A}_G^H$, where $H$ is the function associated to $\overline{\Sigma}_d^\mathfrak{B}$, has infinitely many critical points. 
		
		Next, suppose that $(\overline{\Sigma}_d^\mathfrak{B},\Lambda_d^{\mathfrak{B},E},G)$ is degenerate. In this case, we take a sequence of  Hamiltonians $H_n$   converging to $H$  such that the associated systems are nondegenerate. Then, by the same reason as above, we have
		\begin{equation}\label{RFHn}
			RFH_i^{\mathbb{Z}_2}\big(S^1\times\mathbb{R}^3, H_n^{-1}(0) , L_{0}, G, [\text{pt.}]\big)\cong \mathbb{Z}_2.
		\end{equation}
		We can define spectral invariants for such infinitely many homology classes. By the continuity property of spectral invariants, it follows that $\mathcal{A}^H_G$ also possesses infinitely many critical points, see \cite[Section 4.3]{frauenfelder2019a} for details.
		
		Therefore, we have infinitely many Reeb chords on $\overline{\Sigma}_d^\mathfrak{B}$ joining two Legendrian knots in \eqref{eq:leg}. This yields infinitely many Reeb chords on the original $\overline{\Sigma}_d^\mathfrak{B}$ from $\Lambda_d^{\mathfrak{B},E}$ to $\Gamma_d^{\mathfrak{B},E}$, and in turn, infinitely many symmetric Reeb chords from $\Lambda_d^{\mathfrak{B},E}$ to itself, see Remark \ref{rem:double}. Since such a symmetric Reeb chord corresponds to a symmetric consecutive collision orbit with the Earth, the corollary follows.
	\end{proof}
	
	One weakness of this statement is that infinitely many symmetric Reeb chords obtained in the above proof may unfortunately arise as iterates of a single symmetric Reeb chord that is also periodic. Therefore we can only infer the existence of infinitely many symmetric consecutive collision orbits in the absence of periodic symmetric collision orbits. In the next section,  we improve the above result by showing that certain periodic symmetric collision orbits do not exist under a suitable genericity assumption.

	\subsection{The Existence of at least two Geometrically Distinct Symmetric Consecutive Collision Orbits}

	The rotating Kepler problem (i.e.~$\mu=0$ in \eqref{eq:unreg_H}) can be regularized in the same manner as the R3BP. Since the Hamiltonian of the Kepler problem and the Hamiltonian generating the rotation of the frame (i.e.~the angular momentum) Poisson-commute, every solution of the rotating Kepler problem is simply a rotating image of a solution of the Kepler problem. Therefore, a periodic solution passing through the Earth in the regularized Kepler problem can only occur if its period is a rational multiple of the rotation period of the system. The period $T$ of a periodic solution in the Kepler problem at the energy level $d\in\mathbb{R}$ is given by
	\begin{align*}
		T=\sqrt{\frac{\pi}{2d^3}}=2\pi\sqrt{\frac{1}{8\pi d^3}},
	\end{align*}   
	where the rotation period of the system is $2\pi$. Thus, for almost all energy levels $d\in\mathbb{R}$, no periodic solutions passing through the Earth exist. Adapting the proof of \cite[Theorem~5.1.2]{ruck24b}, we show that this nonexistence also holds for symmetric orbits in the R3BP for a generic mass ratio. 
	
	In the following proposition, in order to reflect the dependence on the mass ratio, we denote by $\Lambda_{d,\mu}^{\mathfrak{B},E}\subset \overline{\Sigma}^\mathfrak{B}_{d,\mu}$ the energy hypersurface and the Legendrian knot associated with mass ratio $\mu\in[0,1]$ at the energy level $d$. 
	
	\begin{figure}[h]
		\begin{center}
			\begin{tikzpicture}
				\draw[arrows=->] (-1.5,0) --(1.5,0);
				\draw[arrows=->] (0,-1.5) --(0,1.5);
				\draw[blue,->] (0, 0) to[out=270, in=315] (-0.71, -0.71);
				\draw[blue] (-0.71, -0.71) to[out=135, in=180] (0, 0);
				\draw[blue,->] (0, 0) to[out=180, in=225] (-0.71, 0.71);
				\draw[blue] (-0.71, 0.71) to[out=45, in=90] (0, 0);
				\draw[blue] (0, 0) to[out=0, in=45] (1.065, -1.065);
				\draw[blue] (1.065, -1.065) to[out=225, in=270] (0.5, 0);
				\draw[blue] (0, 0) to[out=0, in=315] (1.065, 1.065);
				\draw[blue] (1.065, 1.065) to[out=135, in=90] (0.5, 0);
				\draw[white] (4,0) node {$\bullet$};
			\end{tikzpicture}
			\begin{tikzpicture}
				\draw[arrows=->] (-1.5,0) --(1.5,0);
				\draw[arrows=->] (0,-1.5) --(0,1.5);
				\draw[red,->] (0, 0) to[out=45, in=90] (1,0);
				\draw[red] (1, 0) to[out=270, in=315] (0, 0);
				\draw[red] (0, 0) to[out=135, in=180] (0, 1);
				\draw[red] (0,1) to[out=0, in=45] (0, 0);
				\draw[red,->] (0, 0) to[out=225, in=270] (-1, 0);
				\draw[red] (-1, 0) to[out=90, in=135] (0, 0);
				\draw[red] (0, 0) to[out=315, in=0] (0, -1);
				\draw[red] (0, -1) to[out=180, in=225] (0, 0);
			\end{tikzpicture}
		\end{center}
		\caption{Sketch of symmetric periodic orbits with three collisions (left) and four collisions (right).}
		\label{even-odd-collision}
	\end{figure}
	\begin{proposition}
		Let $\epsilon>0$ be as in Theorem \ref{thm:moser2}. Let $d\in(H(L_1),H(L_1)+\epsilon)$ be such that the regularized rotating Kepler problem does not admit periodic solutions passing through the Earth. Then, the set $S\subset (0,1)$ of $\mu$ for which the energy hypersurface $\overline{\Sigma}^\mathfrak{B}_{d,\mu}$ admits a symmetric periodic Reeb orbit intersecting $\Lambda_{d,\mu}^{\mathfrak{B},E}$ an odd number of times is nowhere dense, i.e. for a generic $\mu$ these kind of orbits do not appear.
		\label{noP}
	\end{proposition}
	\begin{proof}
		We begin by examining the constraints that arise from the condition that an orbit is a symmetric orbit with an odd number of intersections with the collision locus. 		
		
		Let $\gamma(\mu,t)$ be such a symmetric periodic Reeb orbit on $\overline{\Sigma}_{d,\mu}^{\mathfrak{B}}$ intersecting  $\Lambda_{d,\mu}^{\mathfrak{B},E}$ an odd number of times. We write 		
		\begin{align*}
			\gamma(\mu,t)=\big(z_1(\mu,t),z_2(\mu,t),w_1(\mu,t),w_2(\mu,t)\big).
		\end{align*}
		In position space, the projection $(z_1,z_2)$ of $\gamma$ goes through the Earth collision point at least once in such a way that it intersects the $z_1$-axis orthogonally. Otherwise, it would intersect the Earth an even number of times. This follows from the fact that, at the Earth collision point, the reflection about the $z_1$-axis does not map any incoming path to the corresponding  outgoing path  and thus collision points appear in pairs due to symmetry. Note that if the reflection does map an incoming path to the corresponding outgoing path, then the orbit  intersects the $z_1$-axis orthogonally at the Earth collision point. See Figure~\ref{even-odd-collision} for a visualisation of orbits with an odd versus an even number of collisions. Therefore,	we may assume that $\gamma(\mu,t)$ is parametrized so that $\gamma(\mu,0)$ intersects the Earth orthgonally. Then, we have 
		\[
		z_1(\mu,0)=\frac{1}{2},\quad z_2(\mu,0)=0,\quad w_1(\mu,0)=0,\quad  w_2(\mu,0)= \pm 2\sqrt{2\mu}.
		\]
		Note that $w_1(\mu,0)=0$ due to the fact that the associated tangent line is symmetric, and $w_2(\mu,0)= \pm 2\sqrt{2\mu}$ follows from $\gamma(\mu,0)\in K^{-1}(0)$, where $K$ is the Hamiltonian in Section \ref{IntroBirk}. Note that two possiblities of the sign of $w_2(\mu,0)$ reflect the fact that both orientations of the $y$-axis are invariant under the symmetry action. We may assume that $w_2(\mu,0)= 2\sqrt{2\mu}$.
		
		Let $S\subset(0,1)$ be the set of masses $\mu$ such that there is  a symmetric periodic Reeb orbit $\gamma(\mu,t)$ on $\overline{\Sigma}_{d,\mu}^{\mathfrak{B}}$ intersecting  $\Lambda_{d,\mu}^{\mathfrak{B},E}$ an odd number of times. Assume by contradiction that there is a sequence with distinct elements $(\mu_n)_{n\in\mathbb{N}}$ in $S$ which converges to some $\mu_*\in S$. Let $\tau_n$ and $\tau_*>0$ be the periods of the orbits $\gamma(\mu_n,t)$ and $\gamma(\mu_*,t)$, respectively. The Hamiltonian equation for $K$ implies that
		\begin{align*}
			\del_t z_2(\mu_*,\tau_*)=\frac{1}{4}w_2(\mu_*,\tau_*) = \frac{\sqrt{2\mu_*}}{2}\neq 0.
		\end{align*}
		Therefore, by the implicit function theorem, there exists  a unique smooth function 
		\[
		\tau:(\mu_*-\delta,\mu_*+\delta)\longrightarrow(0,\infty)
		\] 
		for some $\delta>0$ such that $\tau(\mu_*)=\tau_*$ and 
		\begin{equation}\label{eq:z_2}
			z_2(\mu,\tau(\mu))=0\qquad \forall\mu\in (\mu_*-\delta,\mu_*+\delta).
		\end{equation} 
		Moreover, for sufficiently large $n$, we also have $\tau(\mu_n)=\tau_n$, and therefore,
		\begin{equation}\label{eq:mu_n}
			z_1(\mu_n,\tau(\mu_n))=\frac{1}{2},\qquad   w_1(\mu_n,\tau(\mu_n))=0.	
		\end{equation}
		We consider the functions 
		\[
		f,g: (\mu_*-\delta_0,\mu_*+\delta_0)\longrightarrow \mathbb{R},\qquad f(\mu):=z_1(\mu,\tau(\mu)),\quad  g(\mu):=w_1(\mu,\tau(\mu)).
		\]
		Differentiating $z_2(\mu,\tau(\mu))=0$ with respect to $\mu$, we obtain 
		\begin{equation}\label{eq:d_mu}
			\frac{\D}{\D \mu}\tau(\mu)= -\frac{\del_{1}z_2(\mu,\tau(\mu))}{\del_{2}z_2(\mu,\tau(\mu))}		
		\end{equation}
		where $\partial_1$ and $\partial_2$ denote the partial derivatives with respect to the first and  second coordinates, respectively. The Hamiltonian equation for the restricted three-body problem is real analytic. Since $z_2$ is real analytic, the Cauchy-Kovalevskaya theorem  applied to \eqref{eq:d_mu} implies that $\tau$ is real analytic, and in turn $f$ and $g$ are also real analytic. Since $f-\frac{1}{2}$ and $g$ are zero at $\mu_n$ by \eqref{eq:mu_n}, due to the identity theorem,  they are identically zero. This together with \eqref{eq:z_2} proves that, for every $\mu\in(\mu_*-\delta,\mu_*+\delta)$, $\gamma(\mu,t)$ is a symmetric periodic Reeb orbit on $\overline{\Sigma}_{d,\mu}^{\mathfrak{B}}$ intersecting  $\Lambda_{d,\mu}^{\mathfrak{B},E}$. 
		
		Since $\overline{\Sigma}_{d,\mu}^{\mathfrak{B}}$ is of contact type for all $\mu$, the blue sky catastrophe does not occur for $\gamma(\mu,t)$ smoothly parametrized by $\mu\in(\mu_*-\delta,\mu_*+\delta)$, see \cite[Theorem~7.6.1]{frauenfelder2018a}. Therefore, the periods of $\gamma(\mu,t)$ are uniformly bounded, and as $\mu\to \mu_*-\delta$, $\gamma(\mu,t)$ converges to a symmetric periodic Reeb orbit  $\gamma(\mu_*-\delta,t)$ on $\overline{\Sigma}_{d,\mu_*-\delta}^{\mathfrak{B}}$ intersecting  $\Lambda_{d,\mu_*-\delta}^{\mathfrak{B},E}$.  
		Now we apply the preceding argument, based on the implicit function theorem, to $\gamma(\mu_*-\delta,t)$. By repeating this process, we can successively extend the interval $[\mu_*-\delta,\mu_*+\delta)$ all the way to $\mu=0$. This implies the existence of a symmetric periodic orbit intersecting the collision fiber over the Earth in the case  $\mu=0$, i.e.~in the regularized rotating Kepler problem. However, this contradicts our choice of energy level $d$. This completes the proof.	\end{proof}
	
	\begin{corollary}\label{cor3}
		Let $\epsilon>0$ be as in Theorem \ref{thm:moser2}. Let $d\in (H(L_1),H(L_1)+\epsilon)$ and $\mu\in(0,1)$ be generic in the sense of Proposition \ref{noP}. Then, $\left(\Sigma_{d,\mu},X_{H}\right)$ has at least two (geometrically distinct) symmetric consecutive collision orbits with the Earth that together form a periodic orbit. The same statement holds for the Moon.
	\end{corollary}
	\begin{proof}
		As shown in the proof of Corollary \ref{cor2}, there are 
		infinitely many symmetric Reeb chords on $\overline{\Sigma}_{d,\mu}^\mathfrak{B}$ with endpoints in  $\Lambda_{d,\mu}^{\mathfrak{B},E}$. If none of them is periodic, then they are all geometrically distinct. If this is not the case, due to Proposition \ref{noP}, there is a periodic symmetric collision orbit that intersects the Earth at least twice. Since such an orbit intersects the Earth at least twice, it contains at least two geometrically distinct symmetric consecutive collision orbits. This completes the proof.      
	\end{proof}
	
	\begin{remark}
		Note that the seemingly stronger result \cite[Corollary~5.1.3]{ruck24b} which asserts the existence of infinitely many symmetric consecutive collision orbits for energies below the first critical energy value is also only valid under the additional assumption that there are generically no symmetric periodic orbits that intersect the collision locus an even number of times.
	\end{remark}

	\bibliographystyle{alpha}
	\bibliography{Kbib}	
	
	\end{document}